\documentclass[11pt]{article}
\usepackage[english]{style}
\usepackage{latexsym,amsmath,amssymb,amsfonts,diagrams,graphicx}

\author{%
Andrei C\u ald\u araru, Junwu Tu\thanks{Mathematics Department,
University of Wisconsin--Madison, 480 Lincoln Drive, Madison, WI
53706--1388, USA, {\em e-mail: }{\tt andreic@math.wisc.edu,
tu@math.wisc.edu}}}
\title{Curved $A_\infty$ algebras and Landau-Ginzburg models}
\date{}
\DeclareFontFamily{U}{rsf}{}
\DeclareFontShape{U}{rsf}{m}{n}{
  <5> <6> rsfs5 <7> <8> <9> rsfs7 <10-> rsfs10}{}
\DeclareMathAlphabet{\mathscr}{U}{rsf}{m}{n}

\DeclareMathAlphabet{\mathgth}{U}{euf}{m}{n}

\DeclareFontFamily{U}{cyr}{}
\DeclareFontShape{U}{cyr}{m}{n}{
  <5> wncyr5 <6> wncyr6 <7> wncyr7 <8> wncyr8 <9> wncyr9 <10-> wncyr10}{}
\DeclareMathAlphabet{\mathcyr}{U}{cyr}{m}{n}

\input cyracc.def

\makeatletter
\def\operator@font{\sf}
\makeatother

\newcommand{\cB}{{\mathscr B}}

\newcommand{\cC}{{\mathscr C}}
\newcommand{\cD}{{\mathscr D}}
\newcommand{\cE}{{\mathscr E}}
\newcommand{\cF}{{\mathscr F}}

\newcommand{\cL}{{\mathscr L}}

\newcommand{\cO}{{\mathscr O}}

\newcommand{\dgcat}{{\mathsf{dgcat}}}

\newcommand{\field}[1]{\mathbb{#1}}

\newcommand{\perf}{{\mathsf{perf}}}

\newcommand{\D}{{\mathsf D}}

\newcommand{\tr}{{\mathsf{tr}}}
\newcommand{\BM}{{\mathsf{BM}}}
\newcommand{\DG}{{\mathsf{DG}}}

\newcommand{\tors}{{\mathsf{tors}}}

\newcommand{\bbk}{\mathbf{k}}

\DeclareMathOperator{\Cone}{Cone}
\DeclareMathOperator{\Spec}{Spec}

\DeclareMathOperator{\Der}{Der}

\DeclareMathOperator{\Hom}{Hom}

\DeclareMathOperator{\Coker}{Coker}

\DeclareMathOperator{\Jac}{Jac}
\DeclareMathOperator{\Proj}{Proj}

\DeclareMathOperator{\Tw}{Tw}

\DeclareMathOperator{\imag}{Im}
\DeclareMathOperator{\tot}{tot}
\DeclareMathOperator{\id}{id}

\DeclareMathOperator{\Ext}{Ext}

\newcommand{\ra}{\rightarrow}

\newcommand{\contract}{\lrcorner}

\newcommand{\C}{\mathbb{C}}

\newcommand{\Z}{\mathbb{Z}}

\newcommand{\iso}{\cong}

\newcommand{\pj}{\mathbb{P}}

\newcommand{\HH}{H\mskip-2mu H}

\newcommand{\Mod}{\mathsf{Mod}}
\newcommand{\sg}{\mathsf{sg}}
\newcommand{\Gr}{{\mathsf{Gr}}}
\newcommand{\QGr}{{\mathsf{QGr}}}
\newcommand{\gr}{{\mathsf{gr}}}
\newcommand{\qgr}{{\mathsf{qgr}}}
\newcommand{\Hqe}{{\mathsf{Hqe}}}
\renewcommand{\phi}{\varphi}

\newarrow{Equal} =====

\begin{document}
\maketitle
\begin{abstract}
  We study the Hochschild homology and cohomology of curved $A_\infty$
  algebras that arise in the study of Landau-Ginzburg (LG) models in
  physics.  We show that the ordinary Hochschild homology and
  cohomology of these algebras vanish.  To correct this we introduce
  modified versions of these theories, Borel-Moore Hochschild homology
  and compactly supported Hochschild cohomology. For LG models the new
  invariants yield the answer predicted by physics, shifts of the
  Jacobian ring.

  We also study the relationship between graded LG models and the
  geometry of hypersurfaces.  We prove that Orlov's derived
  equivalence descends from an equivalence at the differential graded
  level, so in particular the CY/LG correspondence is a
  dg equivalence.  This leads us to study the equivariant Hochschild
  homology of orbifold LG models.  The results we get can be seen as
  noncommutative analogues of the Lefschetz hyperplane and Griffiths
  transversality theorems.
\end{abstract} 
\section{Introduction}
\label{sec:intro}

\paragraph
The main goal of this paper is to study the Hochschild homology and
cohomology of categories of matrix factorizations.  This is achieved
by exploiting the relationship that exists between matrix
factorizations and curved associative algebras.  As we shall see, the
usual Hochschild invariants are trivial for curved algebras.  To
correct this we introduce versions of these invariants that are better
behaved in this context.

\paragraph
Categories of matrix factorizations have appeared in several areas of
mathematics and physics.  They were first introduced by Eisenbud in
the study of hypersurface singularities~\cite{Eis}, and further
studied by Buchweitz as the stable category of Cohen-Macaulay
$B/W$-modules~\cite{Buc}.  Two decades later Kontsevich argued that
matrix factorizations are the correct mathematical concept to
formalize the Landau-Ginzburg (LG) models from physics of string
theory.  Subsequently Orlov~\cite{Orlov2} argued that there is a
strong relationship between a graded version of the category of matrix
factorizations and the geometry of projective hypersurfaces.

\paragraph
The main construction in these theories associates to the data of an
associative algebra $B$ and an element $W$ (the superpotential) in the
center of $B$ a $\Z/2\Z$- or $\Z$-graded differential graded (dg)
category $\DG(B, W)$ called the category of matrix factorizations.
Its homotopy category $\D^b(B, W)$ is called the derived category of
matrix factorizations of the pair $(B,W)$.

There is a natural way to associate to the data $(B,W)$ a
$\Z/2\Z$-graded curved $A_\infty$ algebra $B_W$.  As a super-vector
space it consists of $B$, concentrated purely in the even part.  The
only non-zero operations are $m_0$, determined by the equality $m_0(1)
= W$, and $m_2$, given by the original multiplication of $B$.  We will
call such an algebra a curved associative algebra.  

The relationship between the algebra $B_W$ and the category of matrix
factorizations is well understood: the category $\Tw(B_W)$ of twisted
complexes over $B_W$ (the natural generalization to the $A_\infty$
setting of the dg category of complexes of projective modules over an
associative algebra) is precisely $\DG(B, W)$.

From a geometric point of view the $A_\infty$ algebra $B_W$ can be
thought of as the algebra of ``functions'' on a hypothetical
non-commutative space $X$, with $\Tw(B_W)$ playing the role of the
dg category of complexes of coherent sheaves on $X$.

\paragraph
We are primarily interested in understanding the Hochschild homology
and cohomology of the dg category of matrix factorizations $\DG(B, W)
= \Tw(B_W)$.  When dealing with the category of complexes of
projective modules over a ring $A$ (or, more generally, over a
non-curved $A_\infty$ algebra) the situation is greatly simplified by
the existence of the identities
\begin{align*}
\HH^*(\Tw(A)) = \HH^*(A),\\
\HH_*(\Tw(A)) = \HH_*(A).
\end{align*}
Indeed, it is often much easier to compute the invariants of $A$ than
computing those of $\Tw(A)$.

However, as we shall soon see, the above isomorphisms fail for curved
algebras like $B_W$.  In fact the first result we prove in this paper
is the following theorem.

\begin{Theorem} 
\label{thm:van-hh}
For any curved associative algebra $B_W$ we have
\[ \HH^*(B_W) = \HH_*(B_W) = 0. \]
\end{Theorem}
\vspace*{-6mm}

\paragraph
This result should be regarded in a negative way, in the sense that
ordinary Hochschild homology and cohomology are the wrong invariants
to use for curved algebras.  To obtain interesting invariants, we
introduce in Section~\ref{sec:hh} modified versions of these
invariants.  Based on a close analogy with algebraic topology we call
these invariants Borel-Moore Hochschild homology, $\HH^\BM_*$, and
compactly supported Hochschild cohomology, $\HH_c^*$.  They are
closely related to the ordinary Hochschild invariants, but with the
twist that a Borel-Moore chain is an element of the {\em direct
  product} of homogeneous Hochschild chains, as opposed to a direct
sum in the usual definition.  Similarly, a compactly-supported
Hochschild cochain is a direct sum of homogeneous cochains, as opposed
to a direct product.  Collectively we refer to $\HH^\BM_*(B_W)$ and
$\HH_c^*(B_W)$ as the compact type invariants of the curved algebra $B_W$.

\paragraph
The main advantage of the new invariants lies in the fact that an
argument of Segal~\cite{Seg} shows that the compact type invariants of
any algebra $A$, curved or not, agree with the compact type invariants
of $\Tw(A)$.  As before this greatly simplifies the computation of the
compact type invariants for categories of the type $\Tw(A)$. In
particular for the category of matrix factorizations $\DG(B_W) =
\Tw(B_W)$ we get the following theorem.

\begin{Theorem}
\label{thm:jaccomp}
Let $Y=\Spec B$ be a smooth affine scheme of dimension $n$ and assume
that $W\in B$ is a regular function on $Y$ with isolated critical
points.  Then we have 
\begin{align*}
& \HH^*_c(\DG(B_W)) = \HH^*_c(B_W) = \Jac(W)[0], \\
& \HH_*^\BM(\DG(B_W)) = HH_*^\BM(B_W) = \omega(W)[n],
\end{align*}
where $\Jac(W)$ denotes the Jacobi ring of $W$, $\omega(W)$ denotes
the canonical module for $\Jac(W)$, and for a vector space $M$ the
notation $M[i]$ represents the $\Z/2\Z$ graded vector space which is $M$
in degree $i \bmod 2$ and 0 otherwise.
\end{Theorem}

\paragraph
\textbf{Remark.} This result agrees precisely with computations of
the closed string sector of Landau-Ginzburg models in physics.

\paragraph
There are two drawbacks to the new invariants.  First, as in the
topological situation, the compact type Hochschild homology and
cohomology are not homotopy invariant.  The more serious issue is
that in general we do not understand the relationship between the
compact type invariants and the traditional ones.  However, it turns
out that for many interesting $A_\infty$ algebras or categories the
compact type invariants agree with the usual ones.  We call such
algebras or categories of {\em compact type}.
\medskip

\noindent
The easiest examples of compact type algebras are given by the
following theorem.

\begin{Theorem}
  Let $A$ be a finite-dimensional, $\Z$-graded $A_\infty$ algebra
  supported in non-positive degrees.  Then $A$ is of
  compact type. In other words we have
\[ HH^*(A) \iso HH^*_c(A),\quad HH_*(A) \iso HH_*^\BM(A). \]
\end{Theorem}
\vspace*{-4mm}

\paragraph
As we have seen from the explicit calculations in
Theorems~\ref{thm:van-hh} and~\ref{thm:jaccomp}, $B_W$ is not of
compact type.  On the other hand, for a complete regular local ring
$B$, Dyckerhoff~\cite{Dyc} computes directly the ordinary Hochschild
(co)homology of the category $\DG(B, W)$, and his answers agree with
our computations of compact type invariants for the curved algebras
$B_W$.  Putting these results together gives the following corollary.

\begin{Corollary}
\label{cor:mfcpttype}
The dg category $\DG(B, W)$ associated to a commutative, complete
regular local ring $B$ and to a non-unit $W\in B$ with isolated
critical points is of compact type.  We have
\[ \HH_*(\DG(B_W)) \iso \HH_*^\BM(B_W) \iso \HH_*^\BM(\DG(B_W)). \]
\end{Corollary}
\vspace*{-6mm}

\paragraph
The above fact is based on calculations that are quite indirect.  A
direct proof of Corollary~\ref{cor:mfcpttype} using curved Koszul
duality will appear in a separate paper~\cite{Tu} of the second
author, thus giving an alternative approach to the calculations of
Dyckerhoff.

\paragraph
\textbf{The graded case.}  The categories of matrix factorizations
described above share many characteristics with {\em affine} geometry.
The natural generalization to the projective setting is realized by
the so-called categories of {\em graded} matrix factorizations.  To be
precise, assume that the ring $B$ that appears in the definition of
matrix factorizations is $\Z$-graded, and that the superpotential $W$
is homogeneous.  Then we shall say that we are dealing with a graded
LG-model, and there is an appropriate modification $\DG_\Z(B, W)$ of
the definition of $\DG(B, W)$ which takes the grading into account.
The resulting category is called the dg category of graded matrix
factorizations, with associated homotopy category $\D^b_\Z(B, W)$.  It
can be viewed as a particular case of the graded category of
singularities introduced by Orlov~\cite{Orlov2}.  
\medskip

\noindent
The main result of~[loc.cit.] is the following theorem.  We state it
using the original notation, which will be reviewed in
Section~\ref{sec:dg}.  The main objects it relates are two
triangulated categories associated to a Gorenstein graded ring $A$,
$\D^\gr_\sg(A)$ and $\D^b(\qgr\mbox{-}A)$.  Their relationship to
$\D^b_\Z(B,W)$ and $\D^b(Proj B/W)$ will be discussed
in~(\ref{subsec:Orlgrmf}) below.

\paragraph
\textbf{Theorem}
{\it
{\em (Orlov~\cite{Orlov2}).} 
\label{thm:orlov}
Let $A$ be a graded Gorenstein algebra with Gorenstein parameter
$a$.  Then for any integer $i\in\Z$ we have:
\begin{itemize}
\item[(A)] if $a>0$ there exists a semi-orthogonal decomposition
\[ \D^b(\qgr\mbox{-}A)=\langle A(-i-a+1),...,A(-i), \D_\sg^\gr(A)\rangle;\]

\item[(B)] if $a<0$ there exists a semi-orthogonal decomposition
\[ \D_{\sg}^{\gr}(A)=\langle k(-i),... k(-i+a+1),\D^b(\qgr\mbox{-}A)\rangle \]

\item[(C)] if $a=0$ there is an equivalence 
\[\D_{\sg}^{\gr}(A)\cong \D^b(\qgr\mbox{-}A).\]
\end{itemize}
}

\paragraph
\label{subsec:Orlgrmf}
To connect the above theorem to our context, assume that $B$ is a
regular graded ring, $W$ is a homogeneous element of $B$, and denote
by $A$ the graded Gorenstein ring $B/W$.  Orlov's construction yields
a graded category of singularities $\D_\sg^\gr(A)$, and
in~\cite{Orlov2} Orlov shows that this triangulated category agrees
with the derived category of graded matrix factorizations
$\D^b_\Z(B_W)$.

Moreover, if $A$ is commutative (so that we can talk about $\Proj A$),
the category $\D^b(\qgr\mbox{-}A)$ agrees with the derived category of
coherent sheaves on $\Proj A$.  Hence the above theorem ultimately
relates the derived category of graded matrix factorizations
$\D^b_\Z(B, W)$ with the derived category of coherent sheaves on the
hypersurface $\Proj (B/W)$.

\paragraph
Since we are interested in Hochschild invariants, and these do not
behave well with respect to derived equivalence, we need to extend
Orlov's result~(\ref{thm:orlov}) to the dg setting. The result we get
is the following theorem. (The notion of semi-orthogonal decomposition
for dg categories is defined in Appendix~\ref{app:dg}.)

\begin{Theorem}
\label{thm:dg-equiv}
Let $A$ be a Gorenstein algebra with Gorenstein parameter $a$. Then
for any integer $i\in\Z$ we have:
\begin{itemize}
\item[(A)] if $a>0$ there exists a semi-orthogonal decomposition
\[ D_i=\langle\pi A(-i-a+1),...,\pi A(-i), T_i\rangle;\]

\item[(B)] if $a<0$ there exists a semi-orthogonal decomposition
\[ T_i=\langle q k(-i),...q k(-i+a+1),D_i\rangle; \]

\item[(C)]\ if $a=0$ there is an equivalence 
\[ T_i\cong D_i.\]
\end{itemize}
Here $D_i$ and $T_i$ are dg enhancements for $\D^b(\qgr\mbox{-}A)$ and
$\D_\sg^\gr(A)$ respectively, and $\pi$ and $q$ are analogues of the
functors with the same name introduced by Orlov.
\end{Theorem}

\paragraph
The main idea of the proof is to use the notion of dg quotient (due to
Drinfel'd~\cite{Drin}) and its universal properties.  (After we proved
our theorem, we became aware of the main result of Lunts and Orlov's
paper~\cite{LO} which provides an alternative, more general way of
proving the same theorem.)

\paragraph
The next step is to understand graded matrix factorizations from the
point of view of non-commutative geometry.  Our goal is to get a
description of $\DG_\Z(B, W)$ similar to the equivalence
\[ \DG(B, W) \iso \Tw(B_W) \] 
we had in the ungraded case.  The idea is to exploit the observation
from physics~\cite{Wal} that graded LG-models can be realized by an
orbifold construction from ungraded ones.  From a mathematical
standpoint this statement is best expressed as the following result.

\begin{Theorem}
\label{thm:mftwist}
  Let $(B,W)$ be a $\Z$-graded Landau-Ginzburg model with $\deg W =
  d$. Then there exists a $\Z$-graded, curved $A_\infty$ category
  $B_W\sharp\Z_d$ such that the category $\Tw(B_W\sharp\Z_d)$ of
  associated twisted complexes is dg equivalent to $\DG_\Z(B, W)$.
\end{Theorem}

\paragraph 
We are mainly interested in understanding the Hochschild homology of
the category $\DG_\Z(B, W)$, and as in the ungraded case we will do
this by understanding instead the Borel-Moore homology of
$B_W\sharp\Z_d$.  This is justified by a graded analogue of
Corollary~\ref{cor:mfcpttype}, which will also appear in the
forthcoming paper~\cite{Tu}.  The corresponding result states that in
the situation we are interested in (hypersurfaces in projective space)
there is an isomorphism of graded vector spaces
\[ HH_*(\DG_\Z(B,W)) \iso HH_*^\BM(B_W\sharp\Z_d). \]

\paragraph
We state the result of our calculation of $\HH_*^\BM(B_W\sharp\Z_d)$ in
a more general form as a localization formula. The notations in the
following theorem are explained in Section~\ref{sec:ehh}.

\begin{Theorem}
\label{thm:hhlgo}
Let $G$ be a finite group acting on a smooth affine scheme $Y=\Spec
B$, and let $W$ be a $G$-invariant global function on $Y$. Denote by
$[(Y,W)/G]$ the dg category of $G$-equivariant matrix factorizations
of $(B,W)$.  Then we have
\[ HH_*^\BM([(Y,W)/G])\iso \left (\bigoplus_{g \in G}
  HH_*^\BM(Y^g,W|_{Y^g})\right)_G \] 
where for $g\in G$, $Y^g$ is the $g$-invariant subspace of $Y$ and the
subscript $G$ denotes taking coinvariants of the induced $G$ action.
\end{Theorem}

\paragraph
Combining the above calculation with the dg equivalences of
Theorem~\ref{thm:dg-equiv} gives a new proof of Griffiths' theorem on the cohomology of a
hypersurface $X = \Proj B/W$ of projective space $\Proj B$.  Indeed,
since Hochschild homology is invariant under dg equivalences, we have
\[ \HH_*(X) \iso \HH_*(\DG_\Z(B,W)) \iso \HH_*^\BM(B_W\sharp\Z_d), \]
and these latter groups are computed by Theorem~\ref{thm:hhlgo}.  On
the other hand, the Hochschild homology groups of $X$ are
traditionally computed using the well-known
Hochschild-Kostant-Rosenberg isomorphism~\cite{Swa}, and the answer
depends directly on the computation of the Hodge numbers $h^{p,q}(X)$.
These are classically computed using the Lefschetz hyperplane and
Griffiths transversality theorems.

We interpret the conclusion of the previous paragraph as giving an
interpretation of Theorems~\ref{thm:hhlgo} and~\ref{thm:dg-equiv} as 
non-commutative analogues of the Lefschetz hyperplane and Griffiths
transversality theorems from classical complex geometry. 

\paragraph
The paper is organized as follows. In Section~\ref{sec:twist}, for the
purpose of fixing notation, we recall the basics of $A_\infty$
algebras along with details of the twist construction and its
relationship to categories of matrix factorizations.  We also prove
Theorem~\ref{thm:mftwist}.  The usual definitions of Hochschild
homology and cohomology are presented in Section~\ref{sec:hh} in a
way which extends to the curved case as well.  We then prove
Lemma~\ref{lem:van} which states that any $A_\infty$ algebra with only
$m_0$ term has vanishing Hochschild (co)homology.  This is the basic
ingredient in the proof of Theorem~\ref{thm:van-hh}.

In Section~\ref{sec:lg} we specialize to the case of the curved
associative algebras that appear in the study of affine LG models. We
explicitly calculate both the usual Hochschild (co)homology groups and
their compact type variations. The usual ones are shown to vanish and the
modified ones yield the answers predicted by physics.

Section~\ref{sec:dg} is devoted to extending Orlov's results on
derived categories of graded singularities to the dg setting.  We begin by
reviewing Orlov's theorem that relates the derived category of
singularities to the derived category of coherent sheaves on the
corresponding projective variety. We also relate the dg versions of the 
category of graded singularities and the category of graded matrix factorizations.

Motivated by Theorem~\ref{thm:mftwist} we calculate in
Section~\ref{sec:ehh} the Borel-Moore Hochschild homology for LG
orbifolds of the form $B_W\sharp G$.  We prove a localization formula
for homology in this general situation.  Unlike the rest of the paper
where homology and cohomology are on equal footing, our localization
result only applies to homology.  We leave the calculation of
compactly supported Hochschild cohomology of orbifolds for future
work.

For the reader's convenience we collect in Appendix~\ref{app:dg}
certain results on dg quotients of dg categories that we use in the
course of this work. We also define the notion of semi-orthogonal
decomposition of dg categories and prove Lemma~\ref{lem:adm} that is
used in the proofs in Section~\ref{sec:dg}.

\paragraph
\textbf{Acknowledgments.} The authors are indebted to Tony Pantev and
Tobias Dyckerhoff for many useful insights and helpful conversations.
This material is based upon work supported by the National Science
Foundation under Grant No.\ DMS-0901224.

\section{A-infinity categories and the twist construction}
\label{sec:twist}

In this section we recall the notion of $A_\infty$ structure and
explain the twist construction of $A_\infty$ algebras (or, more
generally, $A_\infty$ categories) following the exposition of
Seidel~\cite{Seid}.  Using the twist construction it is easy to see
that the category of matrix factorizations can be identified with the
twist of a curved algebra. The situation for the graded case is more
complicated and we explain how to obtain the category of graded
matrix factorizations by an orbifold construction. The conclusion is
that only the finite group ($\Z/d\Z$) is necessary, as opposed to
$\C^*$ as one might have expected for a $Z$-grading.

We will work over a ground field $\bbk$ and we make once and for all
the choice to work with either $\Z$- or $\Z/2\Z$-graded vector spaces.
All tensor products are understood to be in the symmetric monoidal
category of these graded vector spaces, with the Koszul convention for
signs. If $V$ and $W$ are graded vector spaces, $\Hom(V,W)$ is the
graded vector space of all linear maps $V\ra W$, not necessarily
homogeneous. A map $f\in \Hom(V, W)$ is said to have degree $i$ if it
is homogeneous of degree $i$, i.e., it is a homogeneous degree zero
map $V \ra W[i]$.

\paragraph
Let $A$ be a vector space. Form the tensor coalgebra
\[ B(A):=T^c(A[1])= \bigoplus_{i=0}^{\infty} A[1]^{\otimes i} \] 
endowed with the coproduct
\[ \Delta(a_1|\cdots|a_i) = \Sigma_{j=0}^{i}
(a_1|\cdots|a_j)\otimes(a_{j+1}|\cdots|a_i), \] 
where if the index $j$ is $0$ or $i$ we use the unit of the ground
field as the empty tensor product.

\begin{Definition} 
  An $A_\infty$ algebra structure on $A$ consists of a degree one
  coderivation $m : B(A) \rightarrow B(A)$ of $B(A)$ such that $m\circ
  m=0$.
\end{Definition}

\paragraph
\label{subsec:defm}
As the coalgebra $B(A)$ is cogenerated by $A[1]$ giving a
linear map from $B(A)$ to $A$ is equivalent to giving a coderivation
of $B(A)$.  In other words $C^*(A)$ can be identified with the space
of coderivations on $B(A)$.  

Let $m$ be an $A_\infty$ algebra structure on $A$. Under the above
identification, it corresponds to a homogeneous degree one cochain in
$A$ which we shall denote by $m$ as well. The component of $m$ of
tensor degree $i$ is denoted by $m_i$.  It can be regarded as a map
\[ m_i : A^{\otimes i} \ra A[2-i]. \]

We remind the reader the standard terminology used in the theory of
$A_\infty$ algebras. The term $m_0$ is usually assumed to be zero in
most definitions of $A_\infty$ algebras.  Algebras for which $m_0 = 0$
will be called flat, while those for which $m_0$ does not vanish will
be said to be curved.  We emphasize that most of the existing
homological algebra constructions for flat $A_\infty$ algebras do not
generalize in an obvious way to the curved case.

For flat algebras the degree one map $m_1$ is called the differential
as $m_1^2=0$ in this case.  The map $m_2$ is usually called the
product of the algebra and the terms $m_i$ for $i\geq 3$ are referred
to as higher multiplications.

\paragraph
The equation $m\circ m=0$ translates into an infinite system of
quadratic relations between the $m_i$'s that are sometimes taken as
the definition of $A_\infty$ algebras. The first few relations are
\begin{align*}
& m_{1}\circ m_{0}(1) =0,\\ 
& m_{1}\circ m_{1}(-) \pm m_{2}(m_{0}(1),-) \pm m_{2}(-,m_{0}(1)) =0,\\ 
& m_{1}\circ m_{2}(-,-) \pm m_{2}(m_{1}(-),-) \pm m_{2}(-,m_{1}(-))\\
& \pm m_{3}(m_{0}(1),-,-) \pm m_{3}(-,m_{0}(1),-) \pm
m_{3}(-,-,m_{0}(1)) =0,\\ 
&\vdots 
\end{align*}
where the signs are determined by the Koszul sign convention.

\paragraph
For the purposes of this paper we shall be primarily interested in
curved associative algebras or categories. These are special cases of
$A_\infty$ algebras where the only nonzero multiplications are $m_0$
and $m_2$. To keep with Orlov's notation denote the vector space of
such an algebra by $B$.  The $A_\infty$ relations require $m_2$ to
define an ordinary associative algebra structure on $B$ and the image
of $m_0$ is a one dimensional subspace of $B$ generated by an element
$W$ in the center of $B$. Conversely, the data of an associative
algebra $B$ and an element $W$ in its center determine a curved
algebra which will be denote by $B_W$.  (Usually we shall assume that
the entire algebra is concentrated in the even part of the
$\Z/2\Z$-graded vector space $B$.)

\paragraph
Let $\cC$ be an $A_\infty$ category, i.e., a category with higher
composition operations such that the $A_\infty$ identities hold
whenever the compositions make sense. One can think of
$A_\infty$ algebras as $A_\infty$ categories with one object.

The twist construction of an $A_\infty$ category $\cC$ produces
another $A_\infty$ category $\Tw(\cC)$ which generalizes the notion of
the category of perfect complexes associated to an ordinary algebra.
We recall the construction as explained in~\cite{Seid}. 

\paragraph
By definition the morphism sets of $\cC$ are graded and this grading
is central to the twist construction. We shall denote by $H$ the
grading group, which can be $\Z/2\Z$ or $\Z$.  The construction of
$\Tw(\cC)$ proceeds in two steps.  
\medskip

\noindent\textsl{Step 1.} Form an additive enlargement $\Sigma\cC$ of
$\cC$, which adds direct sums and grading shifts of objects of $\cC$.
\medskip

\noindent Objects of $\Sigma \mathcal{C}$ are formal finite sums of the form
\[ X=\bigoplus_{f\in F} X_{f}[\sigma_{f}]\]
for a finite index set $F\subset \Z$ and objects $X_{f}$ of $\cC$, $\sigma_{f}
\in H$. The morphisms in $\Sigma\cC$ are defined by
\[\Hom_{\Sigma \cC}(\oplus_{f\in F} X_f[\sigma_f],\oplus_{g\in G} Y_g[\tau_g]) =\bigoplus_{f,g}\Hom_\cC(X_{f},Y_{g})[\tau_{g}-\sigma_{f}]\]
where $[-]$ denotes shifting the $H$-grading.  Composition of maps in
$\Sigma\cC$ is defined using matrix multiplication with signs.
\medskip

\noindent\textsl{Step 2.} Form the category $\Tw_H(\cC)$.
\medskip

\noindent 
The objects in $\Tw_H(\cC)$ are twisted complexes, which are defined
to be pairs $(X,\delta)$ consisting of an object $X\in Ob(\Sigma
\cC)$ and a morphism $\delta \in Hom(X,X)$ of degree one such that the
generalized Maurer-Cartan equation holds: 
\[\bigoplus_{i\geq 0} m_{i}(\delta,\cdots,\delta)=0.\]
The morphisms in $\Tw_H(\cC)$ are the same as in $\Sigma \cC$, but
compositions involve insertions of $\delta$.

\paragraph
\textbf{Remark.}  In the second step there is usually an upper
triangular property for $\delta$ which ensures that the Maurer-Cartan
equation makes sense (finiteness of summation).  We will not worry about this issue because the
$A_{\infty}$ categories we will be considering have finitely many
non-zero $m_k$'s. Secondly, one usually takes the idempotent
completion of our definition of $\Tw_H(\cC)$ to have better homological
properties. For the applications in this paper, this completion is not
important.

\paragraph
As an easy example we consider the twist construction applied to a
curved associative algebra $B_W$ viewed as a curved $A_\infty$
category with one object, graded by $H=\Z/2\Z$ with all of $B$
concentrated in the even part.

First recall the classical definition of $\DG(B_W)$, the dg category
of matrix factorizations of $W$. Objects of this category are pairs
$(E,Q)$ where $E$ is a free, $\Z/2\Z$-graded $B$-module of finite rank,
and $Q$ is an odd $B$-linear map on $E$ such that $Q^2=W\cdot\id$. The
set of morphisms between two objects $(E,Q)$ and $(F,P)$ is simply
$\Hom_B(E,F)$, the space of all $B$-linear maps.

There is a differential $d$ on $\Hom_B(E,F)$ which makes $\DG(B_W)$
into a differential graded category. It is defined by the formula
\[ d(\phi):= P\circ \phi -(-1)^{|\phi|}\phi\circ Q.\]
The identity $d^2=0$ follows from the fact that $W$ is in the center of $B$.

\paragraph
We can identify $\DG(B_W)$ with $\Tw_{\Z/2\Z}(B_W)$, as follows.  The
objects of the additive enlargement $\Sigma B_W$ are $\Z/2\Z$-graded
free $B$-modules of finite rank. The twisted complexes add the odd map
$\delta$, and the Maurer-Cartan equation reduces to 
\[ W\cdot\id+ \delta\circ\delta=0.\]
As the degree of $\delta$ is odd, the above equation defines a matrix
factorizations of $W$ (up to a sign). It is easy to check that
composition of morphisms is the same as in $\DG(B_W)$, so we
conclude that 
\[ \Tw_{\Z/2\Z}(B_W)\iso \DG(B_W).\]

\paragraph
\label{subsec:defgrmf}
Of interest to us is also the category of graded matrix factorizations
$\DG_\Z(B_W)$ introduced by Orlov~\cite{Orlov2}. This category is a
$\Z$-graded dg category associated to a graded algebra $B$ and a
homogeneous potential $W$.  An object of $\DG_\Z(B_W)$ consists of two
graded free $B$-module $E_0$, $E_1$ along with homogeneous maps
\begin{align*} E_0&\stackrel{P_0}{\rightarrow} E_1 \mbox{ and }\\
E_1&\stackrel{P_1}{\rightarrow} E_0
\end{align*} 
of degrees $0$ and $d$, respectively, satisfying the matrix
factorization identity 
\[ P_1\circ P_0 =W\cdot \id, \quad P_0\circ P_1 =W\cdot \id. \] It is
convenient to denote such data by $(E,P)$ with $E=E_0\oplus E_1$ and
$P$ the odd map which satisfies the matrix factorization identity
$P^2=W\cdot\id$. This notation is thus the same as that in the
un-graded case, but we should keep in mind the degree requirements for
$E$ and $P$.

The morphism space in the graded case takes into account the
$\Z$-grading of $B$. Explicitly the morphism space between two objects
$(E,P)$ and $(F,Q)$ is
\[ \bigoplus_k \Hom_{\gr\mbox{-}B}\left (E,F(kd)\right )\] 
where the subscript $\gr\mbox{-}B$ means that we only take homogeneous
$B$-linear maps of degree zero. This space can be characterized as the
$G$-invariant subspace of the space of all $B$-linear maps
$\Hom_B(E,F)$. The differential on the $\Hom$ space is defined by the
same commutator formula as in the un-graded case. 
\medskip

\noindent
Our next goal is to realize $\DG_\Z(B_W)$ as the twist construction of
a curved category. 

\begin{Theorem} 
\label{thm:gradedmf}
Let $B$ be a $Z$-graded algebra, and let $W$ be a
homogeneous element in the center of $B$. Denote the degree of $W$
by $d$. Then there exists a $\Z$-graded, curved $A_\infty$ category
$B_W\sharp\Z/d\Z$ with finitely many objects such that the twisted
category $\Tw_\Z(B_W\sharp\Z/d\Z)$ is isomorphic to the category of
graded matrix factorizations $\DG_\Z(B_W)$. 
\end{Theorem}

\paragraph
\label{subsec:crossprod}
If $G$ is a group acting on an associative algebra $B$ then we can
form the cross product algebra $B\sharp G$.  As a $B$-module it is the
same as $B\otimes_\Z \Z G$, and we shall denote an element $b \otimes
g$ by $b\sharp g$.  The multiplication in $B\sharp G$ is given by the
formula
\[ (a\sharp g)(b\sharp h) =ab^g\sharp gh,\]
where $b^g$ is the result of the action of $g$ on $b$.

If $W$ is a central element of $B$ which is invariant under the action
of $G$ then $W\sharp 1$ is central in $B\sharp G$, so we can form the
curved algebra $(B\sharp G)_{W\sharp 1}$.  For simplicity we shall
denote this algebra by $B_W\sharp G$.

The category of twisted modules over $B_W\sharp G$ can be identified
with the category of $G$-equivariant twisted modules over $B_W$, so we
can think of the resulting theory as an orbifolding of the original
Landau-Ginzburg model $B_W$.

\paragraph
We now return to the context of Theorem~\ref{thm:gradedmf}.  As the
vector space $B$ is $Z$-graded it carries a natural
$\C^*$-action. The cyclic subgroup $\Z/d\Z=\left\{\bar{i}|0\leq i\leq
  d-1 \right\}$ embeds into $\C^*$ by $\bar{i} \mapsto \zeta^i$ where
$\zeta =e^{\frac{2\pi\sqrt{-1}}{d}}$. Thus $\Z/d\Z$ also acts on $B$
by
\[ \bar{i} \rightarrow ( f \mapsto \zeta^{i|f|} f ). \] 
In the following arguments we shall denote the group $\Z/d\Z$ by $G$.

Using the construction of~(\ref{subsec:crossprod}) we produce a
$\Z/2\Z$-graded curved algebra $B_W\sharp G$.  For the proof of
Theorem~\ref{thm:gradedmf} we shall construct a $\Z$-graded curved
category whose total space of morphisms is $B_W\sharp G$.  The new
category shall also be denoted by $B_W\sharp G$, where no risk of
confusion exists.

\paragraph
The idea is to consider $B\sharp G$ as a category with $d$ objects
instead of just an algebra. The objects of this category are the
characters of $G$ (one dimensional representations of $G$) and the
morphisms between two objects $\chi_i$ and $\chi_j$ are given by the
invariant part of the $G$ action on $B$ twisted by $\chi_i$ and
$\chi_j$. Explicitly, in our case, we denote the objects of this
category by
\[(\frac{0}{d}),(\frac{1}{d}),\cdots,(\frac{d-1}{d}).\] 
The morphisms between two objects $(\frac{i}{d})$ and $(\frac{j}{d})$
consist of elements of $B$ of degree $(j-i)\bmod d$. Composition of
morphisms is multiplication in $B$. 

\paragraph
This category has a $\Z$-grading induced from that of $B$. However,
with this grading we can not add the curvature element $W$, since it
is of degree $d$ and not $2$ as needed for a $\Z$-graded
$A_\infty$ category.

To fix this we define a new $\Z$-grading on the category $B\sharp G$.
The new grading of an element $f \in Hom((\frac{i}{d}),(\frac{j}{d}))$
is given by 
\[\hat{f}:= \frac{2(|f|-j+i)}{d}. \]
The result is an integer since we required that $|f|=(j-i)\mod d$. In
fact, $\hat{f}$ is always in $2\Z$ for any morphism $f$.

The curvature element $W$ has degree two in the new grading, and hence
we can define a $\Z$-graded $A_\infty$ category structure $B_W\sharp G$ on
$B\sharp G$ by adding the curvature term $m_0(1) = W$ in the
endomorphism space of every object in $B\sharp G$.  This is the desired
category for Theorem~\ref{thm:gradedmf}

\paragraph
Since the category $\cC = B_W\sharp G$ is $\Z$-graded, we can perform
the $\Z$-graded twist construction to it. Our goal is to show that
$\Tw_\Z(\cC)$ can be identified with the category of graded
matrix factorizations $\DG_\Z(B_W)$.

\paragraph 
\label{subsec:maptwdg}
To identify $\DG_\Z(B_W)$ with $\Tw_\Z(B_W\sharp G)$ we first identify
the objects.  Recall that objects of $\Tw_\Z(B_W\sharp G)$ are pairs
of the form 
\[ \left ( \bigoplus_i \frac{a_i}{d}[k_i], \delta \right) \]
with $a_i, k_i \in \Z$, and $\delta$ an endomorphism of hat-degree one
satisfying the Maurer-Cartan equation.  Define a map on the level
of objects by the formula
\[ \left(\oplus_i \frac{a_i}{d}[k_i],\delta \right ) \mapsto (E_0\oplus E_1, P),\] 
where
\begin{align*} 
E_0&:= \bigoplus_i B\left(\frac{k_i d}{2}+a_i\right),\\ 
E_1&:= \bigoplus_i B\left(\frac{(k_i+1)d}{2}+a_i\right), \mbox{ and }\\ 
P &:= \mbox{the corresponding matrix defined by } \delta.
\end{align*}
Here the first summation is over those indices $i$ such that the
corresponding $k_i$ is an even integer, while the second one is over
those for which $k_i$ is odd.

\paragraph
The matrix factorization identity is precisely the Maurer-Cartan
equation. It remains to check that each component of $P$ has the
correct degree. Components of $P$ arise from components of
$\delta$, which are morphisms in $B_W\sharp G$ of hat-degree 1.  In other
words such components are maps in $B_W\sharp G$ of the form
\[ \left(\frac{a}{d}[k]\right) \ra \left (\frac{b}{d}[l]\right), \] 
i.e., elements $x\in B$ whose degree $|x|$ satisfy
\[ l-k+1= \frac{2(|x|-b+a)}{d}.\] 
Note that since the right hand side is necessarily even, $l$ and $k$
are of different parity, and hence the only non-zero components of
the map $P$ will be from $E_0$ to $E_1$ and vice-versa (no
self-maps). Moreover, the above equality can be rewritten as
\begin{align*} |x|&=(\frac{(l+1)d}{2}+j)-(\frac{kd}{2}+i)\\ 
&= (\frac{ld}{2}+j)-(\frac{(k+1)d}{2}+i)+d.
\end{align*} 
The first equality can be applied to the case when $k$ is even and $l$
is odd, and it implies that the degree of the corresponding component
of the resulting map $P_0:E_0\rightarrow E_1$ is zero.  The second
equality similarly shows that if $k$ is odd and $l$ is even, $x$
corresponds to an entry of degree $d$ of the resulting homogeneous
matrix $P_1$ from $E_1$ to $E_0$.

\paragraph 
To check that the morphism sets are the same between the two
categories $\Tw_\Z(B_W\sharp G)$ and $\DG_\Z(B_W)$ we need to check
that every component of a morphism in $\Tw_\Z(B_W\sharp G)$
corresponds to an element of $B$ of (ordinary) degree divisible by
$d$.  Consider such a component, which is a morphism
\[ \phi\in \Hom_{\Tw(B_W\sharp G)}\left(\frac{i}{d}[k],\frac{j}{d}[l]\right)
= \Hom_{B_W\sharp G}\left(\frac{i}{d},\frac{j}{d}\right)[l-k]. \]
We can express its ordinary degree as an element in $B$ using the hat
degree, the natural degree of $B_W\sharp G$:
\[ \hat{\phi}+l-k=\frac{2(|\phi|-j+i)}{d}.\] 
Solving for $|\phi|$ yields
\begin{align*} 
|\phi| &= \frac{(\hat{\phi}+l-k)d}{2}+j-i\\
&=(\frac{(l+1)d}{2}+j)-(\frac{kd}{2}+i)+\frac{\hat{\phi}-1}{2}d
\hspace{.5in} &\mbox{ or }\\
&=(\frac{ld}{2}+j)-(\frac{(k+1)d}{2}+i)+\frac{\hat{\phi}+1}{2}d
\hspace{.5in}&\mbox{ or }\\ 
&=(\frac{ld}{2}+j)-(\frac{kd}{2}+i)
+\frac{\hat{\phi}}{2}d \hspace{.5in}&\mbox{ or }\\
&=(\frac{(l+1)d}{2}+j)-(\frac{(k+1)d}{2}+i)+\frac{\hat{\phi}}{2}d.
\end{align*} 

On the other hand $\phi$ corresponds to a morphism $x$ between shifts
of $B$, as explained in~(\ref{subsec:maptwdg}).  Combining these
facts we see that the degree $|x|$ of $x$ as a graded map in
$\DG_\Z(B_W)$ is given by
\[ \frac{\hat{\phi}-1}{2}d,\  \frac{\hat{\phi}+1}{2}d, \mbox{ or
}\frac{\hat{\phi}}{2}d, \]
according to the parity of $k$ and $l$.

Observe also that the fractions $\frac{\hat{\phi}-1}{2}$,
$\frac{\hat{\phi}+1}{2}$ and $\frac{\hat{\phi}}{2}$ which multiply $d$
are in fact all integers. The reason is that the hat grading is in
$2\Z$ for even shifts, and in $2\Z+1$ after odd shifts. So the
components of every morphism in $\Tw_{\Z}(B_W\sharp G)$ are all
elements of $B$ of degree divisible by $d$, which is precisely the
same as in $\DG_\Z(B_W)$.

Thus we have identified the morphism spaces between the two categories
$\DG_\Z(B_W)$ and $\Tw_\Z(B_W\sharp G)$. It is easy to see that
composition in both categories is given by matrix multiplications and
the differential on $\Hom$ spaces is given by commutators. This
concludes the proof of our identification of the two categories.
\qed

\section{Hochschild-type invariants of curved A-infinity algebras}
\label{sec:hh}

In this section we review the definitions of Hochschild
invariants of $A_\infty$ algebras, ensuring that these definitions
include the curved case.  We present a vanishing result showing
that these invariants vanish for an algebra which only has non-trivial
$m_0$.  Finally we introduce the compact type invariants discussed in
the introduction and give a first example of algebras of compact type.

\paragraph
Let $A$ be an $A_\infty$ algebra.  The space of cochains in $A$ is
the vector space
\[ C^*(A) = \Hom(B(A),A[1]) = \prod_{i=0}^\infty \Hom(A[1]^{\otimes
  i},A[1]),\] 
where $B(A)$ is the free coalgebra generated by $A[1]$.  As
in~(\ref{subsec:defm}) $C^*(A)$ can be identified with the space of
coderivations of $B(A)$.

\paragraph
The space of cochains admits two different gradings: the internal
grading arising from the grading of $A$, and a secondary grading by
tensor degree.  For example, a map $A\otimes A \ra A$ of degree zero
can be regarded as an element of $\Hom(A[1]\otimes A[1], A[1])$ of
internal degree one and of tensor degree two.  In general, a map
$A^{\otimes i} \ra A$ of degree $j$ will have internal degree $i+j-1$.
We shall use the internal grading as the default one.

Note that while the tensor grading is always by integers, the internal
grading is by $\Z$ or $\Z/2\Z$, depending on the grading of $A$.
Therefore the Hochschild (co)homology of $A$ only carries a
$\Z/2\Z$-grading if the algebra $A$ is itself only $\Z/2\Z$-graded (and
it is not a strictly associative algebra, in which case the Hochschild
differentials are homogeneous of degree one with respect to the tensor
grading). 

We emphasize that $C^*(A)$ is a direct {\em product} of $\Hom$-sets,
which is a direct consequence of the fact that $B(A)$ is a direct {\em
  sum}.

\paragraph 
The space $C^*(A)$, being identified with the space of coderivations
on $B(A)$, naturally carries the structure of a Lie algebra with
bracket given by the commutator of coderivations. The $A_\infty$
condition $m\circ m=0$ is equivalent to $[m,m]=0$ as $m$ is of degree
one.

The Jacobi identity implies that the operator $d:C^*(A) \ra C^*(A)$
given by $d=[m,-]$ is a degree one differential ($d^2 = 0$). The
Hochschild cohomology $HH^*(A)$ is defined as the cohomology of the
complex $(C^*(A),d)$.  Note that $d$ is homogeneous with respect to
the internal grading on $C^*(A)$, but not, in general, with respect to
the tensor degree.  The latter happens only if $m_2$ is the only
non-zero structure map. 

\paragraph
Let us denote by $|a|$ the parity of $a$ in the vector space $A[1]$, and let
$\phi$ be a cochain in $\Hom(A[1]^{\otimes k},A[1])$. Then
the explicit formula for the component of $d\phi$ of tensor degree
$l$ is 
\begin{align*}
(d\phi)(a_1|&\cdots|a_l):=\\
= \sum_{j,l\geq k}& (-1)^{j+|a_1|+\cdots+|a_j|}
                        m_{l-k+1}(a_1|\cdots|a_j|\phi(a_{j+1}|\cdots|a_{j+k})|\cdots|a_l)+\\
 + \sum_{i} & (-1)^{|\phi|+|a_1|+\cdots+|a_i|}              
                        \phi(a_1|\cdots|a_i|m_{l-k+1}(a_{i+1}|\cdots|a_{i+l-k+1})|\cdots|a_l).
                        \end{align*}
Note that the second term in the above sum includes contributions from
$m_0$ only if $l=k-1$.  In this case the contribution is given by
\[ \sum_{i=0}^{k-1}(-1)^{|\phi|+|a_1|+\cdots+|a_i|}\phi(a_1|\cdots|a_i|m_0(1)|\cdots a_{k-1}).\]
\paragraph
We are also interested in studying the Hochschild homology of curved
algebras.  The space of chains on $A$ is 
\[ C_*(A):= A \otimes B(A) = \oplus_{i=0}^{\infty} A\otimes A^i. \]
Again, like in the case of the space of cochains, there are two
gradings on $C_*(A)$, an internal one and one given by tensor degree,
and we use the internal grading by default.

In order to define the differential we mimic the standard definition
used in the flat case.  The best way to illustrate this formula is to
use trees drawn on a cylinder as explained in~\cite{KS}.  The formula
one gets for the Hochschild differential
$b:C_*(A) \ra C_*(A)$ is
\begin{align*} 
  b(a_0|\cdots|a_i)= &\sum_{j,k} (-1)^{\epsilon_{k}}
  m_j(a_{i-k+1}|\cdots|a_0|\cdots|a_{j-k-1})|a_{j-k}|\cdots|a_{i-k} +\\
  &\sum_{j,k} (-1)^{\lambda_k}
  a_0|\cdots|a_k|m_j(a_{k+1}|\cdots|a_{k+j})|\cdots|a_i 
\end{align*}
where the signs $\epsilon_k$ and $\lambda_k$ are determined by the Koszul
sign convention: 
\begin{align*} 
  s_i:=& \sum_{0\leq l\leq i} |a_l|,\\
  \epsilon_k =& \sum_{0\leq l \leq k-1} (|a_{n-l}|)(s_i-|a_{n-l}|),\\
  \lambda_k =& \sum_{0\leq l \leq k} (|a_l|). 
\end{align*}
It is straightforward to check that $b^2=0$.  The Hochschild homology
of the algebra $A$ is defined as the cohomology of the complex
$(C_*(A),b)$.

\paragraph
\textbf{Remark.} 
 The key feature of the first summation is that $a_0$ has to
be inserted in $m_j$.  Thus the appearance of $m_0$ does not affect these
terms as $m_0$ allows no insertion.  In the second summation $a_0$ is
not inserted in $m_j$ and has to be placed at the first spot.  Thus if
$m_0$ is present, we should insert the term $m_0(1)$ into any spot
after the first term $a_0$. Explicitly, the terms involving $m_0$ are
\[ \sum_{k=0}^{i} (-1)^{\lambda_k} a_0|\cdots|a_k|m_0(1)|\cdots|a_i.
\]

\paragraph 
From a computational point of view the Hochschild (co)homology of
curved $A_\infty$ algebras is far more complicated than that of
non-curved ones.  The main difficulty is caused by the fact that the
Hochschild (co)homology differential does not preserve the filtration
on (co)chains induced by tensor degree in the curved case.  (It is
easy to see that in the flat case this filtration is preserved, even
though the differentials are not homogeneous with respect to the
tensor degree.)  Thus the usual spectral sequences associated with
this filtration that are used for flat algebras are unusable for
curved ones.  We shall not go into more details as we shall not use
these spectral sequences.


\paragraph
We now present a vanishing result which will motivate the discussion
of compact type invariants for curved algebras later.

\begin{Lemma} 
\label{lem:van} 
Let $A$ be an $A_\infty$ algebra over a field $k$ such that $m_0\neq
0$ and $m_k = 0$ for any $k > 0$.  This is equivalent to the data of a
graded vector space $A$ together with the choice of a degree two
element $W = m_0(1)$ in $A$.  Then both the Hochschild homology and
cohomology of $A$ are zero.
\end{Lemma}

\begin{Proof} 
To avoid sign issues we shall assume that $A$ is concentrated in even
degrees.  The proof still works with corrected signs in the general
case.  Also, throughout this proof the degree of a (co)chain will mean the
tensor degree, so that $C_k(A)$ means $A^{\otimes (k+1)}$.  

We begin with the statement for homology.  The Hochschild chain
complex for our algebra $A$ is
\[ C_*(A) = \bigoplus_i A^{\otimes (i+1)} \] 
with differential given by
\[b(a_0|\cdots|a_i)=\sum_{k=0}^{i} (-1)^{k} a_0|\cdots|a_k|W|
\cdots|a_i. \] 
We will prove that this complex is acyclic by constructing an explicit
homotopy between the identity map and the zero map. 

Let $L:A\ra k$ be a $k$-linear map such that $L(W) \neq 0$.  Such a
map exists as $k$ is a field and $A$ is free as a $k$-module, so we
can extend any nonzero map from the one-dimensional subspace spanned
by $W$ to the whole space $A$.  Once such a functional $L$ is chosen
we can use it to define a homotopy $h_{k}:C_k(A) \ra C_{k-1}(A)$ by
the formula:
\begin{align*} 
h_k(a_0|\cdots|a_k)&=(-1)^{k+1} L(a_k)(a_0|\cdots|a_{k-1}), \\ 
h_0&=0. 
\end{align*}
An explicit calculation shows that $h_*$ is a homotopy between the
identity map on $C_*(A)$ and the zero map, hence $C_*(A)$ is acyclic.



                           

A similar argument yields the result for cohomology, with homotopy
given by $h_{k}:C^k(A)\rightarrow C^{k+1}(A)$,
\begin{equation}
[h_k(\varphi_k)](a_1|\cdots|a_{k+1})=L(a_1)\varphi_k(a_2|\cdots|a_{k+1}). 
\tag*{\qed} 
\end{equation}
\end{Proof}
\medskip

\noindent
Motivated by the vanishing results in the above lemma we
introduce the following modifications of the Hochschild invariants. 

\begin{Definition} 
Let $A$ be an $A_\infty$ algebra or category with finitely many
nonzero operations $m_k$.  Consider modified spaces of chains and cochains 
\begin{align*} 
C_*^\Pi(A) & = \prod_{i=0}^\infty A\otimes A^{\otimes i}, \\ 
C^*_\oplus(A) & = \bigoplus_{i=0}^\infty \Hom(A^{\otimes i}, A). 
\end{align*} 
The Borel-Moore Hochschild homology $HH_*^\BM(A)$ is defined to be the
homology of the complex $(C_*^\Pi(A),b)$, and the compactly
supported Hochschild cohomology $HH_c^*(A)$ is defined to be the
cohomology of $(C_\oplus^*(A),d)$, where $b$ and $d$ are
given by the same formulas as before.
\end{Definition}

\paragraph
\textbf{Remark.} 
Without the assumption that $A$ has finitely many nonzero
higher multiplications these invariants might not be well-defined due
to the possible non-convergence of the infinite sum.  Note that the
property required of $A$ is not homotopy invariant: it is easy to
construct dg algebras ($A_\infty$ algebras that only have non-zero
$m_1$ and $m_2$) whose minimal models have infinitely many non-zero
$m_k$'s.  

\paragraph 
Let $A$ be an $A_\infty$ algebra or category for which the
compact type invariants can be defined.  There exist natural maps 
\[ \HH_c^*(A) \ra \HH^*(A)\quad\mbox{and}\quad \HH_*(A) \ra
\HH_*^\BM(A) \] 
induced by the chain maps $(C_\oplus^*(A), d) \ra (C^*(A), d)$ and
$(C_*(A), b) \ra (C_*^\Pi(A), b)$ given by the natural inclusion of
the direct sum in the direct product.

\begin{Definition}
  When the above maps are isomorphisms we shall say that $A$ is of
  {\em compact type}.
\end{Definition}

\paragraph
It is reasonable to expect that under certain finiteness hypotheses we
get compact type algebras. The following proposition shows
that this is the case in a simple situation.

\begin{Proposition} 
  Let $A$ be a $\Z$-graded, finite dimensional $A_\infty$ algebra
  concentrated in non-positive degrees. Assume that $A$ has only
  finitely many nonzero higher multiplications over a field $k$.  Then
  $A$ is of compact type.  In other words we have
  \begin{align*} 
    HH_*(A)&\cong HH_*^\pi(A), \mbox{ and } \\ 
    HH^*(A)&\cong HH_\oplus^*(A).
\end{align*} 
\end{Proposition}
\vspace*{-6mm}

\begin{Proof}
It suffices to observe that under the assumptions of this proposition, the ordinary Hochschild cochain complex $C^*(A)$ is the completion of the compactly supported Hochschild cochain complex $C^*_c(A)$ with respect to the natural internal $\Z$-grading. The case for Hochschild homology is similar.\qed
\end{Proof}

\section{Hochschild invariants of LG models}
\label{sec:lg}
In this section we compute the usual and compact type Hochschild
invariants of the curved $A_\infty$ algebras that arise in the study
of LG models. More precisely we show that the usual Hochschild
invariants vanish for any curved algebra while the modified
invariants agree with the predictions from physics.

\paragraph
The setup we shall work with is the following.  Let $Y=\Spec B$ be a
smooth affine scheme of dimension $n$ and assume that $W\in B$ is a
regular function on $Y$.  Let $Z$ denote the set of critical points of
$W$, i.e., the subscheme of $Y$ cut out by the section $dW$ of
$\Omega^1_Y$.  The ring of regular functions on $Z$ shall be denoted
by $\Jac(W)$, the Jacobi ring of $W$.  The relative dualizing sheaf
$\omega_{Z/Y}$ of the closed embedding $Z\hookrightarrow Y$ is a
coherent sheaf on $Z$, i.e., a module $\omega(W)$ over $\Jac(W)$.
Algebraically $\Jac(W)$ and $\omega(W)$ are the cokernels of the maps
$\contract dW: \Der(B) \ra B$ and $\wedge dW: \Omega^{n-1}_B \ra
\Omega^n_B$, respectively.

Denote by $B_W$ the curved algebra with $m_2$ given by multiplication in
$B$ and with $m_0$ given by $W$. Since $B$ is commutative $W$ is trivially
in the center and hence $B_W$ is an $A_\infty$ algebra.

\begin{Theorem}
\label{thm:hhHH}
With the above notations, we have
\begin{itemize}
\item[(a)] $HH^*(B_W) = HH_*(B_W) = 0$.
\item[(b)] If $\dim Z = 0$, i.e., the critical points of $W$ are
  isolated, then we have isomorphisms of
  $\Z/2Z$-graded vector spaces
\begin{align*}
HH^*_c(B_W) & = \Jac(W)[0], \\
HH_*^\BM(B_W) & = \omega(W)[n].
\end{align*}
(Here $M[i]$ denotes the $\Z/2\Z$ graded vector space which is $M$ in
degree $i \mod 2$ and 0 otherwise.)
\end{itemize}
\end{Theorem}

\paragraph
\textbf{Remark.}  We see that the usual result that $\HH^*(A) =
\HH^*(\Tw(A))$ for flat (non-curved) $A_\infty$ algebras fails for
curved ones.  Indeed, for the curved algebra $B_W$ we have $\HH^*(B_W)
= 0$, while
\[ \HH^*(\Tw(B_W)) = \HH^*(\DG_{\Z/2\Z}(B_W)) = \Jac(W)[0] =
\HH^*_c(B_W). \]
Here the second equality is due to Dyckerhoff~\cite{Dyc} (but some
more assumptions are needed for $B$ and $W$).  See
also~\cite{Tu}.  The third equality is part (b) of the theorem above.

\paragraph
The differential $b$ on Hochschild chains decomposes as
\[b=b_-+b_+\]
where
\[ b_+:C_n(B)\rightarrow C_{n+1}(B),\ 
b_-:C_{n-1}(B)\rightarrow C_n(B) \]
are given by
\begin{align*}
b_+(a_0|\cdots|a_n) &=\sum_{i=0}^n (-1)^i a_0|\cdots|a_i|W|a_{i+1}|\cdots|a_n, \\
b_-(a_0|\cdots|a_n) &=\sum_{i=0}^{n-1}(-1)^ia_0|\cdots|a_ia_{i+1}|\cdots|a_n+(-1)^na_na_0|a_1|\cdots|a_{n-1}.
\end{align*}
Here we are using the tensor degree for the spaces of chains or
cochains.

\begin{Lemma}
\label{lem:mix}
We have $b_+^2=0$, $b_-^2=0$, $b_+b_-+b_-b_+=0$.
\end{Lemma}

\begin{Proof}
  The fact that $b_+^2=b_-^2 = 0$ is a straightforward calculation.
  Alternatively for the second equality one could observe that $b_-$
  is the Hochschild differential for the algebra $B$ (without $m_0$)
  for which it is well known that $b_-^2 = 0$.  Finally the last
  equation follows from the first two and the fact that
  $b^2=(b_++b_-)^2=0$ as the Hochschild differential of the
  $A_\infty$ algebra $B_W$.
\end{Proof}

\paragraph
From the above lemma we see that the Hochschild chain complex
$C_*(B_W)$ is actually a so-called mixed complex.
Recall~\cite[9.8]{Weibel} that a mixed complex is a $\Z$ graded
complex with two anti-commuting differentials of degrees $1$ and $-1$
respectively. Mixed complexes appear in the theory of cyclic homology
where the two differentials are the Hochschild boundary map and the
Connes cyclic operator. We find this analogy quite useful and we will
calculate the Hochschild (co)homology of $B_W$ using techniques
developed for cyclic homology.  We refer the reader to Loday's
book~\cite{Lod} for a detailed discussion of cyclic homology.  The
idea for computing the homology of a mixed complex $C$ is to associate
a double complex $BC$ to it, which on one hand avoids the problems
caused by the inhomogeneity of $b$ and on the other hand allows us to
use the powerful machinery of spectral sequences. As the double
complexes we are going to study are usually unbounded, it will make a
difference if we take the direct sum or the direct product total
complex.  To distinguish between the two we shall denote the direct
sum and direct product total complex of a double complex $BC$ by
$\tot^\oplus(BC)$ and $\tot^\Pi(BC)$, respectively. For a bounded
double complex $BC$ the direct sum and the direct product total
complex are the same so we will simply denote any of them by
$\tot(BC)$.

\paragraph
Consider the double complex $BC_{i,j}=C_{j-i}(B_W)=B\otimes
B^{\otimes(j-i)}$ associated to the mixed complex
$(C_*(B_W),b_+,b_-)$. It is a bicomplex supported above the diagonal
on the $(i,j)$-plane as illustrated below (here $i$ and $j$ are
horizontal and vertical coordinates respectively):
\[
\begin{matrix}
\mbox{\space}  &  \cdots            &        \mbox{\space}      &       \cdots      &       \mbox{\space}     &      \cdots      &     \mbox{\space}\\
\mbox{\space}  &  \downarrow        &        \mbox{\space}      &       \downarrow  &       \mbox{\space}     &      \downarrow  &     \mbox{\space}\\
\leftarrow     & C_2(B_W)           &        \leftarrow         &       C_1(B_W)    &       \leftarrow        &      C_0(B_W)    &     \leftarrow  0\\
\mbox{\space}  &  \downarrow        &        \mbox{\space}      &       \downarrow  &       \mbox{\space}     &      \downarrow  &     \mbox{\space}\\
\leftarrow     & C_1(B_W)           &        \leftarrow         &       C_0(B_W)    &       \leftarrow        &      0           &    \leftarrow\cdots\\
\mbox{\space}  &  \downarrow        &        \mbox{\space}      &       \downarrow  &       \mbox{\space}     &      \downarrow  &     \mbox{\space}\\
\leftarrow     & C_0(B_W)           &        \leftarrow         &        0          &       \leftarrow        &      0           &   \leftarrow\cdots\\
\mbox{\space}  &  \downarrow        &        \mbox{\space}      &       \downarrow  &       \mbox{\space}     &      \downarrow  &     \mbox{\space}\\
\leftarrow     &  0                 &        \leftarrow         &        0          &       \leftarrow        &      \cdots      &      \mbox{\space}
\end{matrix}
\]
The horizontal differential of $BC$ is given by $b_+$ and the vertical
one is $b_-$.

Note that $BC$ is periodic with respect to shifts along the main
diagonal, hence the two total complexes $\tot^\oplus(BC)$ and
$\tot^\Pi(BC)$ are 2-periodic and this periodicity descends to their
homology.  Therefore it makes sense to talk about the even or the odd
homology of $\tot^\oplus(BC)$ or $\tot^\Pi(BC)$ and for
$*=\mbox{even or odd}$ we have
\begin{align*}
H_{*}( \tot^\oplus(BC) )& \cong HH_{*}(B_W) \mbox{ and }\\
H_{*}( \tot^\Pi(BC))    &\cong HH_*^\BM(B_W).
\end{align*}

\paragraph
We will need another double complex which is the part of $BC$ that
lies in the first quadrant of the $(i,j)$-plane.  Denote this
double complex by $BC^+$. Note that $BC^+$ is bounded, hence the
direct sum and the direct product total complexes agree. The
relationship between $BC^+$ and $BC$ is that the former is a quotient
of the latter by the subcomplex consisting of those terms whose
$i$-coordinate is strictly negative.

\paragraph
\textsl{Proof of part (b) of Theorem~\ref{thm:hhHH}}.  Consider the
positive even shifts $\tot(BC^+)[2r] \mbox{ for } r\in \field{N}$.
For $r<t$ the complex $\tot(BC^+)[2r]$ is a quotient of
$\tot(BC^+)[2t]$ due to the $2$-periodicity.  Hence the quotient maps
\[ \tot(BC^+)[2t]\rightarrow \tot(BC^+)[2r] \]
form an inverse system whose inverse limit is $\tot^\Pi(BC)$
(see~\cite[5.1]{Lod} for more details).  In short we have realized
$\tot^\Pi(BC)$ as the inverse limit
\[ \tot^\Pi(BC)\cong \lim_\leftarrow \tot(BC^+)[2r]. \]

Moreover the tower $(\tot(BC^{+})[2r],r\in \field{N})$ satisfies the
Mittag-Leffler condition~\cite[3.5.6]{Weibel} as these maps are all
onto.  By~\cite[3.5.8]{Weibel} we get a short exact sequence for
$k\in \Z$
\[
0\rightarrow\lim_\leftarrow{}^1 H_{k+1}(\tot(BC^{+})[2r])\rightarrow H_{k}(\tot^{\Pi}(BC))\rightarrow \lim_\leftarrow H_{k}(\tot(BC^{+})[2r])\rightarrow 0.
\]
Since we have
\[ H_*(\tot(BC^+)[2r]) \cong H_{*+2r}(\tot(BC^+)), \]
the problem of computing
\[ HH_*^\BM(B_W) = H_*(\tot^\Pi(BC)) \] 
reduces to the problem of computing the homology of $\tot(BC^{+})$.
Since the bicomplex $BC^{+}$ is non-zero only in the first quadrant,
the standard bounded convergence theorem~\cite[5.6.1]{Weibel} applies.
We will use the spectral sequence associated to the vertical
filtration of the bicomplex $BC^{+}$ whose ${}^1 E$ page consists
precisely of the Hochschild homology of the associative ring $B$.
(The vertical differential $b_-$ is precisely the Hochschild
differential of $B$.)  Since $B$ is regular the classical
Hochschild-Kostant-Rosenberg isomorphism applies and gives
\[ {}^1 E_{ij} = \Omega^{j-i}_B. \]
Moreover, the horizontal differential at the ${}^1 E$ page is induced
from the original horizontal differential $b_+$ on $BC^+$, and the
former can be calculated using the splitting map $e$ of the
HKR-isomorphism, see~\cite[9.4.4]{Weibel}
\[e_{k}(a_0|\cdots|a_k)=\frac{1}{k!}a_0da_1\wedge\cdots\wedge da_k.\]
Since we have
\begin{align*}
e_{k+1}\circ b_+(a_0|\cdots|a_k)&= e_{k+1}(\sum_{i=0}^k (-1)^{i} a_0|\cdots| a_i|W|\cdots|a_k)\\
&=\frac{1}{(k+1)!}\sum_{i=0}^k (-1)^i a_0\wedge\cdots\wedge da_i\wedge dW\wedge \cdots \wedge da_k\\
&=\frac{1}{k!}dW\wedge a_0\wedge da_1\wedge\cdots\wedge da_k
\end{align*}
we conclude that the horizontal differential on the ${}^1 E$ page is the
map $\Omega^{j-i}_B \ra \Omega^{j-i+1}_B$ given by $\alpha \mapsto
dW\wedge \alpha$.


To calculate the ${}^2 E$ page, we compute homology with respect to
this horizontal differential.  Observe that the rows in the ${}^1 E$
page are precisely truncations of the dual of the Koszul complex
associated to the section $dW$ of $\Omega^1_B$.  The assumption that
the critical points of $W$ are isolated implies that this complex is
exact except at the $n$-th (last) spot (recall that $n=\dim B$).

Therefore the ${}^2 E$ page is everywhere zero except for the spots
$(i,i+n)$ for $i>0$ and the spots $(0,j)$ for $0\leq j \leq
n$. More precisely it is $(\wedge^{j}\Omega_{B})/(dW)$ at the spots
$(0,j), 0\leq j \leq \mbox{dim}_k\mbox{ } B$ and $\omega_B/(dW)$ at
the spots $(i,i+\mbox{dim}(B)), \forall i\in \field{N}^{+}$. So the
spectral sequence degenerates at the $E^{2}$ page already.  We
conclude that the inverse limit of the homology groups of
$\tot(BC^{+})[2r]$ is

\[
\lim_\leftarrow H_*(\tot(BC^+)[2r])\cong \begin{cases}
                                              \omega_B/(dW) & *=\mbox{dim}_k\mbox{ } B (\mbox{mod } 2)\\
                                              0             & \mbox {otherwise.}

                                \end{cases}
\]

It also follows that the maps in the tower of these homology groups
are actually isomorphisms for large $r$. Therefore the tower of homology
groups also satisfies the Mittag-Leffler condition, which implies the
vanishing of the first derived functor of the inverse limit
\[
\lim_{\leftarrow}\displaystyle^{1}H_{k}((\tot BC^{+})[2r])=0,\  \forall k\in \Z
\]
Hence we conclude that
\[
HH_*^\BM(B_W)\cong H_*(\tot^\Pi(BC)) \cong  \begin{cases}
                                              \omega_B/(dW) &
                                              *\equiv\dim_k B \bmod 2)\\
                                              0             & \mbox {otherwise.}

                                                 \end{cases}
\]
\paragraph
A similar argument can be adapted for the computation of $HH^*_c(B_W)$.
In fact, in this case, the proof is easier as we will not need the argument 
involving inverse limits. The conclusion is that

\[HH_\oplus^*(B_W)\cong H_*(\tot^{\oplus}(BC)) \cong \begin{cases}
                                              \mbox{Jac}(W) & *=\mbox{even}\\
                                              0      & *=\mbox{odd}.
                                     \end{cases}
                                     \]

\paragraph
\textsl{Proof of part (b) of Theorem~\ref{thm:hhHH}.} The usual Hochschild homology 
of $B_W$ can be computed by the direct sum total complex of $BC$. For this, the convergence 
theorem of spectral sequences shows that the spectral sequence associated to the horizontal 
filtration converges to $\tot^\oplus(BC)$. To calculate the ${}^1 E$ page of this spectral sequence, 
we can apply Lemma~\ref{lem:van}, which shows that the ${}^1 E$ page is already all zero. Thus we have
proved that
\[HH_{*}(B_W) \cong H_{*}(\tot^\oplus(BC)) \cong 0.\]

\paragraph
For the case of the ordinary Hochschild cohomology of $B_W$, 
again by Lemma~\ref{lem:van}, we can deduce that the usual Hochschild 
cohomology vanish by an argument using the idea of taking inverse limit.

The usual Hochschild cohomology $HH^*(B_W)$ was defined to be the
homology of $d=d_++b_-$ on the direct product of Hochschild cochains
$\prod_{n\geq 0} C^n(B_W)$.  Thus we should consider $\tot^{\Pi}(BC)$
of the double complex $BC$ associated to the mixed complex
$(C^*(B_W),d_+,d_-)$. Again the vertical filtration realizes
$\tot^\Pi(BC)$ as the inverse limit of the tower of complexes
$\tot(BC^+)[2r]$ for $r \in \field{N}$. Hence we have the following
short exact sequence:
\[
0\rightarrow\lim_\leftarrow\displaystyle^1 H_{k+1}(\tot(BC^{+})[2r])\rightarrow H_{k}(\tot^{\Pi}(BC))\rightarrow \lim_\leftarrow H_{k}(\tot(BC^{+})[2r])\rightarrow 0
\]
To compute the homology of $\tot(BC^+)$ we run the spectral sequence
associated to the vertical filtration of $BC^+$. Then the first page
of the spectral sequence can be calculated using the vertical
differential $d_-$. But according to Lemma~\ref{lem:van} this is
already zero. Using the above short exact sequence, the vanishing of
the homology of $\tot(BC^+)$ immediately implies the vanishing of
$H_k(\tot^\Pi(BC))$.  We conclude that
\begin{equation}
HH^{*}(B_W)\cong H_*(\tot^\Pi(BC))\cong 0.\tag*{\qed}
\end{equation}

\paragraph
\textbf{Remark.} The above proof of the vanishing of the usual
Hochschild homology and cohomology actually works for any curved
algebra over any field of arbitrary characteristic.

\section{Orlov's results in the dg setting}
\label{sec:dg}

In~\cite{Orlov2} Orlov relates the derived category of graded matrix
factorizations of a homogeneous superpotential $W$ to the derived
category of coherent sheaves on the projective hypersurface defined by
$W = 0$.  In this section we argue that Orlov's results can be
extended to the dg setting by exhibiting a chain of quasi-equivalences
between appropriate dg enhancements of the categories he considered.
The main context in which we work is that of the category $\Hqe$ (see
Appendix~\ref{app:dg}) in which quasi-equivalences are inverted, and
thus in $\Hqe$ the resulting categories are isomorphic.

\paragraph 
We begin by reviewing Orlov's results from~\cite{Orlov2}.  Let
$A=\oplus_{i\geq 0} A_i$ be a connected graded Gorenstein Noetherian
$k$-algebra. Recall the the Gorenstein property means that $A$ has
finite injective dimension $n$, and $R\Hom_A(k,A)$ is isomorphic to
$k(a)[-n]$ for some integer $a$ which is called the {\em Gorenstein
  parameter} of $A$. Here $()$ and $[]$ are the internal degree shift
and homological degree shift respectively.

Denote by $\gr\mbox{-}A$ the abelian category of finitely generated
graded left $A$-modules.  Let $\tors\mbox{-}A$ be the full subcategory
of $\gr\mbox{-}A$ consisting of all graded $A$-modules that are finite
dimensional over $k$ and form the quotient abelian category
$\qgr\mbox{-}A:=\gr\mbox{-}A/\tors\mbox{-}A$.  When $A$ is
commutative $\qgr\mbox{-}A$ is equivalent to the category of coherent
sheaves on $\Proj A$.

For convenience we will sometimes omit the algebra $A$ in our
notation.  For example $\gr\mbox{-}A$ will simply be denoted by
$\gr$.

There is a canonical quotient functor $\pi$ from $\gr$ to $\qgr$,
which in the commutative case corresponds to the operation associating
a sheaf on $\Proj A$ to a graded $A$-module. Moreover $\pi$ admits a
right adjoint $\omega$ defined by
\[ \omega(N):= \oplus_{i\in\Z}\Hom_{\qgr\mbox{-}A}(A, N(i)). \]
Again, in the geometric setting, $\omega$ is the analogue of the
functor $\Gamma_*$ of~\cite[II.5]{HarAG}.

For any integer $i$ we also consider the abelian category $\gr_{\geq
  i}$ which is the full subcategory of $\gr$ consisting of graded
$A$-modules $M$ such that $M_p=0$ for $p<i$.  The categories
$\tors_{\geq i}$ and $\qgr_{\geq i}$ are similarly defined, and the
restrictions of the functors $\pi$ and $\omega$ to these subcategories
shall be denoted by $\pi_{\geq i}$ and $\omega_{\geq i}$.

Deriving these abelian categories yields triangulated categories
$D^b(\gr)$, $D^b(\gr_{\geq i})$ and $D^b(\qgr)$.  However, when
defining the derived functor of $\omega_{\geq i}$ we run into trouble
as the category $\qgr$ has neither enough projectives nor enough
injectives. We can however bypass this difficulty by giving up the
condition of our modules being finitely generated.  Namely consider
the category $\Gr$ which consists of all graded $A$-modules, finitely
generated or not, and form the corresponding quotient $\QGr$. Note
that $\D^+(\qgr)$ is naturally a full subcategory of $\D^+(\QGr)$, and
similarly for $\D^b$, etc.

Since $\QGr$ has enough injectives we can define the derived
functor
\[R\omega_{\geq i}: D^+(\QGr) \rightarrow D^+(\Gr_{\geq i}). \] 
It can be checked that the restriction of $R\omega_{\geq i}$ to the
subcategory $D^+(\qgr)$ lands inside $D^+(\gr_{\geq i})$.  If we
assume further that $A$ has finite injective dimension then the
restriction of $R\omega_{\geq i}$ to $D^b(\qgr)$ lands inside
$D^b(\gr_{\geq i})$. One can also show that the functor $\pi_{\geq i}$
is exact (and hence does not need to be derived) and it is a one-sided
inverse of $R\omega_{\geq i}$: $\pi_{\geq i} \circ R\omega_{\geq
  i}\iso \id$.

The main object of study in this section is the quotient of $D^b(\gr)$
by the full triangulated subcategory $\perf(\gr)$ consisting of
perfect objects.  (Recall that a complex is perfect if it is
quasi-isomorphic to a bounded complex of projectives.) This category
is known as the (graded) category of singularities and is denoted by
$D_{\sg}^{\gr}$.

\paragraph
\textbf{Theorem}
{\it
{\em (Orlov~\cite{Orlov2}).} 
\label{thm:orl}
Let $A$ be a Gorenstein algebra with Gorenstein parameter $a$. Then
the triangulated categories $D_{\sg}^{\gr}$ and $D^b(\qgr)$ are
related as follows:
\begin{itemize}
\item[(A)] if $a>0$ then for each integer $i\in\Z$ there exists a fully faithful
functor $\Phi_i:D_{\sg}^{\gr}\rightarrow D^b(\qgr)$ and a
semi-orthogonal decomposition
\[ D^b(\qgr)=\langle\pi A(-i-a+1),...,\pi A(-i), \Phi_i(D_{\sg}^{\gr})\rangle;\]

\item[(B)] if $a<0$ then for each integer $i\in\Z$ there exists a fully faithful
functor $\Psi_i:D^b(\qgr)\rightarrow D_{\sg}^{\gr}$ and a
semi-orthogonal decomposition
\[ D_{\sg}^{\gr}=\langle q k(-i),...q k(-i+a+1),
\Psi_i(D^b(\qgr))\rangle \] 
where $q:D^b(\gr)\rightarrow D_{\sg}^{\gr}(A)$ is the natural
projection;

\item[(C)] if a=0 then
$D_{\sg}^{\gr}\cong D^b(\qgr)$.
\end{itemize}
}

\paragraph 
Our purpose is to generalize the above result to the dg setting so we
can apply invariance of Hochschild structures under dg equivalences.
For this we need appropriate dg enhancements of the triangulated
categories involved. For the category $D^b(\qgr)$ we take the dg
category of complexes of injective objects in $\QGr$ that are
quasi-isomorphic to objects in $D^b(\qgr)$.  We denote this dg
category by $\DG(\qgr)$.

For the category $D_{\sg}^{\gr}$ there is no obvious candidate for a
dg enhancement. However observe that $D_{\sg}^{\gr}$ is obtained as
the quotient of $D^b(\gr)$ by $\perf(\gr)$, and these two categories
admit natural dg enhancements by using either injective or projective
resolutions of the corresponding complexes. Both dg enhancements will
be used in this paper and the relationship between them will be
clarified in Lemma~\ref{lem:rel}.

We first consider the category obtained by injective resolutions. Let
$\underleftarrow{\DG}(\gr)$ be the dg category of bounded below
complexes of injective graded $A$-modules that are quasi-isomorphic to
objects of $D^b(\gr)$ (the left arrow is used to indicate that we are
using injective resolutions).  Denote by $\underleftarrow{P} (\gr)$
the full subcategory of $\underleftarrow{\DG} (\gr)$ whose objects are
perfect complexes. Then the dg enhancement of $D_{\sg}^{\gr}$ that we
will consider is the dg quotient $\underleftarrow{\DG}
(\gr)/\underleftarrow{P} (\gr)$, which we shall denote by
$\underleftarrow{\DG_\sg^\gr}$. Similarly, we denote by
$\underleftarrow{\DG}(\gr_{\geq i})$ the dg version of $D^b(\gr_{\geq
  i})$ and by $\underleftarrow{P_{\geq i}}(\gr)$ the dg version of
$\perf_{\geq i}(\gr)$.

\begin{Lemma}
\label{lem:inj} The functor $\omega_{\geq i}$ from $\QGr$ to
$\Gr_{\geq i}$ sends injectives to injectives.
\end{Lemma}

\begin{Proof} 
This is immediate from the fact that the left adjoint $\pi_{\geq i}$
of $\omega_{\geq i}$ is exact.
\qed
\end{Proof}

\paragraph 
Because of this lemma we have a well-defined dg functor
\[\omega_{\geq i}: \DG(\qgr)\rightarrow \underleftarrow{\DG}(\gr_{\geq
i}). \]
We denote by $D_i$ the full subcategory of
$\underleftarrow{\DG}(\gr_{\geq i})$ consisting of objects in the image
of $\omega_{\geq i}$.  Note that the functor $\omega_{\geq i}$ induces
a quasi-equivalence between $\DG(\qgr)$ and $D_i$ by the corresponding
results in the triangulated setting. Moreover, the subcategory $D_i$
is an admissible subcategory of $\underleftarrow{\DG}(\gr_{\geq i})$
as the adjoint to the inclusion functor is given by $\pi_{\geq i}$.

\begin{Lemma} 
The dg functor 
\[ \iota_i: \underleftarrow{\DG}(\gr_{\geq i})/\underleftarrow{P_{\geq
    i}}(\gr)\rightarrow \underleftarrow{\DG} (\gr)/\underleftarrow{P}
(\gr)\] 
induced by the natural morphism between localization pairs 
\[ (\underleftarrow{P_{\geq i}}(\gr), \underleftarrow{\DG}(\gr_{\geq
  i})) \ra(\underleftarrow{P} (\gr),\underleftarrow{\DG} (\gr)) \]
is a quasi-equivalence. (See Appendix~\ref{app:dg} for details on
localization pairs).

Moreover, the full subcategory $\underleftarrow{P_{\geq i}}(\gr)$ is
left admissible in $\underleftarrow{\DG}(\gr_{\geq i})$. Thus by
Lemma~\ref{lem:adm} in the Appendix the category $T_i:={}^\bot
\underleftarrow{P_{\geq i}}(\gr)$ in $\underleftarrow{\DG}(\gr_{\geq
i})$ is canonically quasi-isomorphic to the quotient
$\underleftarrow{\DG}(\gr_{\geq i})/\underleftarrow{P_{\geq i}}(\gr)$.
\end{Lemma}

\begin{Proof} 
We only need to show that the morphism $\iota$ induces an equivalence
at the level of homotopy categories. Since all our dg categories are
pre-triangulated, the problem is reduced to that of showing that the
morphism 
\[ D^b(\gr_{\geq i})/\perf_{\geq i} \ra D^b(\gr)/\perf(\gr)\] 
is an equivalence, a result which is proved by Orlov in~\cite{Orlov2}.

The second statement also follows from the definition and from results
of Orlov. 
\qed
\end{Proof}

\paragraph 
We have two admissible full subcategories $D_i$ and $T_i$ inside
$\underleftarrow{\DG}(\gr_{\geq i})$ such that $D_i$ is isomorphic to
$\DG(\qgr)$ and $T_i$ is isomorphic to
$\underleftarrow{\DG_\sg^\gr}(A)$ in $\Hqe$. Once these categories
have been put inside the same dg category Orlov's proof of the
existence of the corresponding semi-orthogonal decompositions
(Theorem~\ref{thm:orl}) also yields a proof of the following
dg analogue, as the notion of orthogonality is strictly on the level
of objects (see Lemma~\ref{lem:adm} in the Appendix~\ref{app:dg}). The
proof of the following theorem will be omitted as it is identical to
Orlov's original proof of Theorem~\ref{thm:orl}.

\begin{Theorem} Let $A$ be a Gorenstein algebra with Gorenstein
parameter $a$. Then we have:
\begin{itemize}
\item[(A)] if $a>0$ there exists a semi-orthogonal decomposition
\[ D_i=\langle\pi A(-i-a+1),...,\pi A(-i), T_i\rangle;\]

\item[(B)] if $a<0$ there exists a semi-orthogonal decomposition
\[ T_i=\langle q k(-i),...q k(-i+a+1),D_i\rangle; \]

\item[(C)] if $a=0$ there exists an equivalence $T_i\cong D_i$.
\end{itemize}
\end{Theorem}
\medskip

\noindent
The following result is a reformulation of the above theorem in the
geometric setting, corresponding to a similar theorem stated by Orlov
in~\cite{Orlov2} for triangulated categories (the LG/CY correspondence).

\begin{Theorem} Let $X$ be a Gorenstein projective variety of
  dimension $n$ and let $\cL$ be a very ample line bundle such that
  $\omega_X \cong \cL^{-r}$ for some integer $r$. Suppose that
  $H^j(X,\cL^k)=0$ for all $k\in \Z$ when $j\neq 0,n$. Set
  $A=\oplus_{i\geq 0} H^0(X,\cL^i)$. Then we have
\begin{itemize}
\item[(A)] if $r>0$ ($X$ is Fano) then there is a semi-orthogonal
  decomposition of dg categories
\[\DG(X)=\langle\cL^{-r+1},\cdots,\cO_X, T\rangle\] 
for some $T$ that is isomorphic to $\underleftarrow{\DG_\sg^\gr}$ in
$\Hqe$.

\item[(B)] if $r<0$ ($X$ is of general type) then there is a semi-orthogonal
decomposition of dg categories
\[\underleftarrow{\DG_\sg^\gr}=\langle k(r+1),\cdots,k,D\rangle\] 
for some $D$ that is isomorphic to $\DG(X)$ in $\Hqe$.

\item[(C)] if $r=0$ ($X$ is Calabi-Yau) then there is an equivalence in $\Hqe$
\[\underleftarrow{\DG_\sg^\gr} \cong \DG(X). \] 
\end{itemize}
\end{Theorem}

\paragraph 
We also want to relate the categories considered above to the category
of graded matrix factorizations from Section~\ref{sec:twist}.  Let
$B=\oplus_{i\geq 0}B_i$ be a finitely generated connected graded
algebra over $k$ and let $W\in B_n$ be a central element of degree $n$
which is not a zero-divisor. Consider the quotient graded algebra
$A:=B/W$, the quotient of $B$ by the two-sided ideal generated by $W$.
We are interested in the category of graded matrix factorizations
$\DG_\Z(B_W)$. 

\paragraph
To relate $\DG_\Z(B_W)$ to the dg category of graded singularities we
will use a different model of the latter, given by projective
resolutions instead of injective ones. Denote by
$\underrightarrow{\DG}(\gr)$ the dg category whose objects are bounded
above complexes of free $A$-modules that are quasi-isomorphic to
objects of $D^b(\gr)$, and let $\underrightarrow{P}(\gr)$ be the full
subcategory consisting of perfect objects. As before we denote by
$\underrightarrow{\DG}_\sg^\gr(A)$ the corresponding quotient.

\paragraph
For any dg category $\cD$, we can associate a dg category $\cD_{\leq 0}$ whose
objects are the same as that of $\cD$. The morphism space $\cD_{\leq 0}$ between
two objects $E$ and $F$ is the smart truncation of the complex $\Hom_\cD(E,F)$
in nonpositive degrees. The inclusion functor
\[ \cD_{\leq 0} \ra \cD\]
is a homotopy equivalence. Next we define a dg functor 
\[\Coker: \DG_\Z(B_W)_{\leq 0} \ra  \underrightarrow{\DG}(\gr).\]
For this purpose we think of objects in $\DG_\Z(B_W)$ as graded periodic complexes
\[E = \cdots E_1\stackrel{p_{1}}{\ra} E_0 \ra E_1(n) \ra E_0(n) \ra \cdots.\]
The functor $\Coker$ is defined by
\[E\mapsto \big( \cdots E_1 \ra E_0 \ra 0\ra \cdots \big)\otimes_{B} A,\] 
that is the object $\Coker(E)$ is the truncation of the graded periodic complex $E$
followed by tensoring with the ring $A=B/W$. The reason this functor
is called $\Coker$ is that $\Coker(E)$ gives a natural
projective resolution over $A$ of the Cohen-Macaulay module $E_0/\imag
p_{1}$.

We can compose the functor $\Coker$ with the natural universal
quotient functor $Q:\underrightarrow{\DG}(\gr)\rightarrow
\underrightarrow{\DG}(\gr)/\underrightarrow{P}(\gr)=\underrightarrow{\DG}_\sg^\gr(A)$
in $\Hqe$ to form a morphism in $\Hqe$
\[ Q\circ \Coker: \DG_\Z(B_W)_{\leq 0} \rightarrow
\underrightarrow{\DG}_\sg^\gr(A).\]
\medskip

\noindent
Again, by passing to homotopy categories, Orlov's results
immediately imply the following theorem, which further implies that
the two dg categories $\DG_\Z(B_W)$ and $\underrightarrow{\DG}_\sg^\gr(A)$
are isomorphic in $\Hqe$.

\begin{Theorem} The morphism $Q\circ \Coker$ is an isomorphism in
$\Hqe$.
\end{Theorem}

\paragraph 
As a final part of this section we want to argue that there is an
isomorphism in $\Hqe$
\[ \underrightarrow{\DG}_\sg^\gr(A) \iso
\underleftarrow{\DG}_\sg^\gr(A). \] 
For this consider the dg category $\DG_\gr(\Mod\mbox{-}A)$ of
complexes of graded $A$-modules that are quasi-isomorphic to objects
in $D^b(\gr)$. Let $P$ be the full subcategory of
$\DG_\gr(\Mod\mbox{-}A)$ consisting of perfect objects.  (Note that by
definition $P$ contains all the acyclic objects.) We can consider the
following two morphisms of localization pairs:
\begin{align*}
j:(\underleftarrow{P}(\gr),\underleftarrow{\DG}(\gr))&\rightarrow
(P,\DG_\gr(\Mod\mbox{-}A)),\\
k:(\underrightarrow{P}(\gr),\underrightarrow{\DG}(\gr))& \rightarrow
(P,\DG_\gr(\Mod\mbox{-}A)).
\end{align*}

\begin{Lemma}
\label{lem:rel} 
Both $j$ and $k$ induce quasi-equivalences on the corresponding
quotient categories. Thus $\underrightarrow{\DG}_\sg^\gr(A)$ and
$\underleftarrow{\DG}_\sg^\gr(A)$ are isomorphic in $\Hqe$.
\end{Lemma}

\begin{Proof} 
Consider the induced map on the homotopy categories of the
quotients. As these categories are all pre-triangulated, we have
\begin{align*} H^0(\underleftarrow{\DG}_\sg^\gr(A))&\cong
H^0(\Tw(\underleftarrow{\DG}_\sg^\gr(A)))\\ &\cong
H^0(\Tw(\underleftarrow{\DG}(\gr)))/H^0(\Tw(\underleftarrow{P}(\gr)))\\
&\cong H^0(\underleftarrow{\DG}(\gr))/H^0(\underleftarrow{P}(\gr))\\
&\cong D_\sg^\gr(A).
\end{align*}
We also have
\[ H^0(\DG_\gr(\Mod\mbox{-}A)/P)\cong D_\sg^\gr(A).\]
Thus we have shown that $j$ induces a quasi-equivalence on the
quotients. The proof of the corresponding statement for $k$ is
similar.
\qed
\end{Proof}

\section{Hochschild homology of LG orbifolds}
\label{sec:ehh}

In this section we give a localization formula for Borel-Moore
Hochschild homology groups of algebras of the form that arise in the
study of orbifold LG models.

\paragraph
\label{subsec:motiv}
The main motivation for this calculation is the following.  Let $X =
\Proj B/W$ be a smooth projective Calabi-Yau hypersurface of degree
$d$ in $\pj^{d-1}= \Proj B$, $B= \bbk[x_0,\ldots, x_{d-1}]$.  By the
results in the previous sections there is a dg equivalence
\[ \DG(X) \iso \DG_\Z(B_W) = \Tw_\Z(B_W\sharp G), \] 
where $\DG(X)$ is a dg enhancement of the derived category of $X$ and
$G = \Z/d\Z$. (See Section~\ref{sec:dg} for the appropriate results in
the non-CY case.) Hochschild homology is invariant under
dg equivalences hence
\[ \HH_*(X) = \HH_*(\Tw_\Z(B_W\sharp G)) \] 
and the latter groups can be related to the Borel-Moore homology groups
$\HH_*^\BM(B_W\sharp G)$ (see~\cite{Tu}), which are the groups we
compute in this section.

\paragraph 
Let $G$ be a finite group acting on a smooth affine scheme $Y=\Spec B$
and let $W$ be a $G$-invariant global function on $Y$. Then we can
form the curved cross product algebra $B_W\sharp
G$~(\ref{subsec:crossprod}).  We shall denote the Borel-Moore
Hochschild homology of this algebra by $\HH_*^\BM([(Y,W)/G])$ to
emphasize the geometric point of view.

\begin{Theorem}
\label{thm:hhlgo1} 
In the above context there exists a graded vector space isomorphism
\[ HH_*^\BM([(Y,W)/G])\cong \left (\bigoplus_{g \in G}
HH_*^\BM(Y^g,W|_{Y^g})\right )_G, \] 
where $Y^g$ is the $g$-invariant subspace of $Y$ and the subscript $G$
means taking coinvariants of the induced $G$
action.
\end{Theorem}

\begin{Proof}
We first recall the following result of Baranovsky~\cite{Bar} for the
non-curved case. Let $Y=\Spec B$ be an affine scheme with a finite
group $G$ acting on it. Then we have
\[ HH_*([Y/G]) \cong (\oplus_{g \in G} HH_*(Y^g))_G. \]
In fact this result is stronger in the sense that there is a natural
chain level map that induces the above isomorphism.  The map $\Psi$ that
goes from the Hochschild chain complex of the cross product algebra
$B\sharp G$ to the $G$-coinvariants of the direct sum of the
Hochschild chain complexes on the invariant subspaces is given explicitly
by
\[\Psi (a_0\sharp g_0|\cdots|a_n\sharp
g_n)=(a_0|g_0(a_1)|g_0g_1(a_2)|\cdots|g_0\cdots g_{n-1}(a_n))_g \]
where the subscript $g=g_0\cdots g_n$ means that all the functions on the
right hand side are viewed as functions restricted on $Y^g$.

To see that $\Psi$ indeed defines a map of chain complexes we note
that the Hochschild differential consists of direct sums of face maps
$d_i$, and $\Psi$ commutes with all these face maps except the last
one ($d_n$). However one can check by a direct computation that
$\Psi$ also commutes with $d_n$ after taking coinvariants.

Theorem~\ref{thm:hhlgo1} can now be proved quite easily with the
quasi-isomorphism $\Psi$ in hand. Let us denote by $B^g$ the global
regular functions on the subspace $Y^g$.  (This notation is slightly
misleading: in fact $B^g$ is the set of {\em coinvariants} in $B$ of the
action of $g$; we keep our notation for geometric purposes.)  Denote
by $B_W^g$ the curved algebra $(B^g)_W$.

Consider the Hochschild chain complexes
\begin{align*} C_*^\Pi(B_W\sharp G) &=\prod_n C_n(B_W\sharp G) \mbox{ and }
\\ \left (\oplus_{g\in G} C_*^\Pi(B_W^g)\right )_G &= \left (\bigoplus_{g\in G} \prod_n
C_n(B_W^g)\right )_G
\end{align*} 
which are both mixed complexes with differentials given by the
Hochschild differential $b_-$ and the differential $b_+$ coming from
inserting $W$.  The map $\Psi$ defines a linear map between the
associated double complexes of these mixed complexes which commutes
with $b_-$.
A direct computation shows that $\Psi$ commutes with $b_+$
as well:
\begin{align*}
b_+&(\Psi(a_0\sharp g_0|\cdots|a_n\sharp g_n)) = b_+(
(a_0|g_0(a_1)|g_0g_1(a_2)|\cdots|g_0\cdots g_{n-1}(a_n))_g ) \\ &=
\sum_{i=0}^n (-1)^i (a_0|g_0(a_1)|\cdots|g_0\cdots
g_{i-1}(a_i)|W|\cdots|g_0\cdots g_{n-1}(a_n))_g.
\end{align*} 

\begin{align*} 
\Psi&(b_+(a_0\sharp g_0|\cdots|a_n\sharp
g_n))=\Psi(\sum_{i=0}^{n}(-1)^i (a_0\sharp g_0|\cdots a_i\sharp
g_i|W\sharp e|\cdots|a_n\sharp g_n)) \\ &= \sum_{i=0}^n (-1)^i
(a_0|g_0(a_1)|\cdots|g_0\cdots g_i(W)|\cdots|g_0\cdots
g_{n-1}(a_n))_g\\ &=\sum_{i=0}^n (-1)^i (a_0|g_0(a_1)|\cdots|g_0\cdots
g_{i-1}(a_i)|W|\cdots|g_0\cdots g_{n-1}(a_n))_g; \\
\end{align*}

We conclude that $\Psi$ is a map of double complexes. Taking the
spectral sequences associated to the vertical filtrations of these
double complexes as in Section~\ref{sec:lg} yields a map between
spectral sequences which we shall denote by $\Psi$ as
well. Baranovsky's isomorphism shows that the induced map on the $E^1$
page is everywhere an isomorphism. It follows from the comparison
theorem of spectral sequences (see for example~\cite{Weibel}) that the
two spectral sequences converge to the same homology groups. Since the
first one converges to $HH_*^\BM([(Y,W)/G])$ by definition, while the
second one converges to $(\oplus_{g \in G} HH_*^\BM(Y^g,W|_{Y^g}))_G$
by the results in Section~\ref{sec:lg}, the result is proved.
\qed
\end{Proof}

\paragraph 
\textbf{Example.}  Consider the case where $B=\bbk[x_1,\ldots, x_d]$,
$Y= \Spec B = \C^d$, and $W\in B$ is a homogeneous polynomial of
degree $d$.  We take $G=\Z/d\Z$ acting diagonally on $Y$.  Assume that
$X = \Proj B/W$ is smooth.  Then Theorem~\ref{thm:hhlgo1} yields
\[ \HH_*^\BM([(Y,W)/G])= \left( \HH_*^\BM(Y,W)\right)_G \oplus \C \oplus \cdots
\oplus \C, \] 
where we have $d-1$ copies of $\C$ indexed by the non-trivial elements
in the group $G$.  (The contribution of the non-trivial elements of
$G$ to orbifold homology is traditionally called the twisted sectors
of the theory.)  Note that by the considerations
in~(\ref{subsec:motiv}) the above is also a computation of $\HH_*(X)$.

The above calculation can be regarded as a non-commutative version of
the Lefschetz hyperplane theorem (for the part involving the twisted
sectors) and of the Griffith transversality theorem (for the
$G$-coinvariant part of the identity component, which can be computed
using a Jacobian ring calculation as in Section~\ref{sec:hh}).

For example for a smooth quartic surface in $\pj^{3}$ its even
Hochschild homology groups are of dimensions $1$, $22$, $1$,
respectively, obtained by adding the even verticals in the Hodge
diamond below:
\[ 
\begin{diagram}[height=1em,width=1em]
& & 1 & & \\
& 0 & & 0 & \\
1 & & 20 & & 1 \\
& 0 & & 0 & \\ 
& & 1. & & 
\end{diagram}
\]
\medskip

\noindent
The degree $0$, $4$, $8$
components of the Jacobian ring of the defining quartic have
dimensions $1$, $19$, $1$ respectively. Compared with the above Hodge
diamond what is missing is another $3$ dimensional vector space down
the middle vertical (which corresponds to $\HH_0$), which are
precisely the $3$ one-dimensional twisted sectors.  (The odd Hochschild
homology is trivial on both sides.)

\appendix
\section{Dg categories and dg quotients}
\label{app:dg}

In this appendix we collect some results on dg categories and on the
notion of dg quotients. Our main reference for dg categories is
Keller's paper~\cite{Kell3}.  For dg quotients we refer to the
original construction of Keller~\cite{Kell2}, see also 
Drinfel'd~\cite{Drin} for an alternative construction.

\paragraph 
We shall fix a ground field $k$ and all categories will be assumed to
be $k$-linear categories. (All the constructions below can be
generalized to any ground ring after suitable flat resolutions.)

A dg category is a category $\cD$ such that the $\Hom$ spaces are
complexes over $k$.  Composition maps are required to be not only
$k$-linear but also maps of complexes. Explicitly, for any objects $X$,
$Y$, $Z\in \cD$ the homomorphism space $\Hom_\cD(X,Y)$ is a
$k$-complex and the composition map
\[\Hom_\cD(X,Y)\otimes \Hom_\cD(Y,Z) \rightarrow \Hom_\cD(X,Z)\] is a
map of complexes. Here the tensor product is the tensor product of
complexes.

\paragraph
Similarly a dg functor $F:\cC \rightarrow \cD$ is defined to be a
$k$-linear functor such that the map
\[ \Hom_\cC(X,Y)\rightarrow \Hom_\cD(FX,FY) \] 
is a map of complexes for any objects $X$ and $Y$ of $\cC$.

\paragraph
One can associate to a dg category $\cD$ its homotopy category
$H^0(\cD)$ whose objects are the same as those of $\cD$ and the homomorphism
set between two objects is defined by
\[\Hom_{H^0(\cD)}(X,Y)=H^0(\Hom_\cD(X,Y)).\] 
Moreover the composition in $H^0(\cD)$ is defined to be the one
induced from $\cD$ (well-defined as the original composition maps are
maps of complexes and hence induce maps on the cohomology). It
is easy to see that dg functors induce maps on the corresponding
homotopy categories. For a dg functor $F$ we shall denote the induced
functor on homotopy categories by $H^0(F)$.

\paragraph
A dg functor $F:\cC\rightarrow \cD$ between dg categories $\cC$ and
$\cD$ is said to be a quasi-equivalence if $H^0(F)$ is an equivalence.
Two dg categories are said to be quasi-equivalent if they belong to
the same equivalence class with respect to the equivalence relation
generated by the above notion of quasi-equivalence.

\paragraph
Consider the category of small dg categories over $k$, denoted by
$\dgcat_k$.  Its objects are small dg categories and its morphisms are
dg functors. As quasi-equivalences induce natural isomorphisms on most
homological algebra constructions we would like to consider a
modification of $\dgcat_k$ wherein quasi-equivalences are inverted.
This is more or less in the same spirit as the construction of derived
categories where quasi-isomorphisms between complexes are inverted.
This construction can be carried out as explained in~\cite{Kell3}, and
the resulting category $\Hqe$ is obtained as the localization of
$\dgcat_k$ with respect to quasi-equivalences.

One can show that various types of homological algebra invariants of
dg categories factor through $\Hqe$. These include for example
Hochschild homology, cyclic homology, Hochschild cohomology and more
sophisticated invariants like open-closed string operations
constructed by Costello~\cite{Co1}.

In particular in Section~\ref{sec:dg} we proved the existence of a
relationship between the dg category of coherent sheaves and the
dg category of matrix factorization, both regarded as objects in
$\Hqe$. In the Calabi-Yau case this implies that these two categories
are isomorphic in $\Hqe$, which further implies that these two
dg categories carry isomorphic homological invariants as mentioned
above. See Section~\ref{sec:ehh} for an application.

\paragraph 
There exists a twist construction for dg categories analogous to the
one for $A_\infty$ categories.  It allows us to relate dg categories
to triangulated categories.  Explicitly the twist of a dg category
$\cD$ is a new dg category $\Tw(\cD)$ such that the homotopy category
of $\Tw(\cD)$ has a natural triangulated structure. We shall denote
this resulting triangulated category by $\cD^\tr$.

For example, an ordinary algebra $A$ can be seen as a category with
only one object. Then $\Tw(A)$ is the dg category of bounded complexes
of free $A$-modules. Its homotopy category is in general not the
derived category of $A$, but rather a fully faithful subcategory
consisting of free objects. This example is a special case of the
general statement that $\cD^\tr$ is the fully faithful subcategory of
$D(\cD)$ consisting of representable objects.

The dg category $\Tw(\cD)$ is also called the pre-triangulated
envelope of $\cD$. Intuitively speaking it is the smallest dg category
containing $\cD$ such that its homotopy category has a triangulated
structure. Thus we shall call a dg category $\cD$ pre-triangulated if
the following two condition hold: 
\begin{itemize} 
\item   For any
  object $X\in \cD$, $X[1]$ is isomorphic to an object of $\cD$
  inside $\cD^\tr$. 

\item For any closed morphism $f:X\rightarrow Y \in \cD$, $\Cone(f)$
  is isomorphic to an object of $\cD$ inside $\cD^\tr$. 
\end{itemize}

\noindent
(The object $X[1]$, if it exists, is the unique object of $\cD$
representing the functor $\Hom_\cD(\,-\,, X)[1]$, where the latter is the
shift-by-one of the complex $\Hom_\cD(\,-\,,X)$. A similar definition
applies to the notion of the cone of a morphism.)

It follow easily from the definition that if $\cD$ is pre-triangulated
then $H^0(\cD)$ has a triangulated structure, and the natural
embedding of $\cD$ into $\Tw(\cD)$ is a quasi-equivalence. The induced
functor on homotopy categories is a triangulated equivalence between
triangulated categories.

\paragraph 
If $\cC$ is a triangulated category and $\cE$ is a full triangulated
subcategory of $\cC$ then the triangulated quotient category $\cC/\cE$
is obtained by localizing the category $\cC$ with respect to the
multiplicative system
\[ S=\left\{f\in \Hom(\cC)~|~ \Cone(f) \in \cE\right\}.\] 
This is known as the triangulated quotient construction. There is a
natural quotient functor from $\cC \ra\cC/\cE$.

There is also a quotient construction in the dg context initiated by
Keller~\cite{Kell2} and later elaborated by Drinfel'd~\cite{Drin}. The
follow theorem summarizes the main results on dg quotients.

\begin{Theorem}
\label{thm:quo} 
Let $\cE$ be a full subcategory of a dg category $\cD$. Then there
exists a dg category $\cD/\cE$ together with a quotient map 
\[ Q:\cD\rightarrow \cD/\cE\] in the category $\Hqe$ such that $Q$ and
$\cD/\cE$ have the universal property that every morphism in $\Hqe$
from $\cD$ to some other dg category $\cF$ that annihilates $\cE$ (the
image of any object in $\cE$ is isomorphic to zero in $H^0(\cF)$)
factors through $\cD/\cE$.
\medskip

\noindent
Moreover, the dg quotient is the dg analogue of the triangulated
quotient: we have
\[ (\cD/\cE)^\tr \cong \cD^\tr/\cE^\tr. \]
\end{Theorem}
\vspace*{-6mm}

\paragraph
Dg quotients have good functorial properties with respect to 
{\em localization pairs}.  A localization pair $\cB$ is a pair
of dg categories $(\cB_1, \cB_2)$ such that $\cB_1$ is a full
subcategory of $\cB_2$. A morphism $F$ between localization pairs
$\cB$ and $\cD$ is a dg functor from $\cB_2$ to $\cD_2$ that sends
objects of $\cB_1$ to objects of $\cD_1$.  Keller showed that such a
map $F$ induces a map on the dg quotients
\[ F: \cB_2/\cB_1 \rightarrow \cD_2/\cD_1.\]

\paragraph
For the remainder of this section we are interested in understanding
semi-orthogonal decompositions at the dg level, and these make best
sense in the case of pre-triangulated dg categories. Let $\cD$ be a
pre-triangulated dg category, and denote by $\cC=H^0(\cD)$ the
homotopy category of $\cD$. By the previous discussion $\cC$ is
triangulated.

Recall the notion of semi-orthogonal decomposition of a
triangulated category $\cC$. Let $\cE$ be a full triangulated
subcategory of $\cC$. We would like to decompose $\cC$ as a ``sum'' of
$\cE$ and its orthogonal complement. In order to do this we need the
subcategory $\cE$ to be {\em admissible}: $\cE$ is said to be right (left)
admissible if the inclusion functor from $\cE$ to $\cC$ admits a
right (left) adjoint.

Assume that $\cE$ is a right admissible full triangulated
subcategory of $\cC$.  Define the right orthogonal complement of
$\cE$ to be the full subcategory $\cE^\bot$ of $\cC$ consisting of those
objects $X$ $\in \cC$ such that
\[ \Ext^*(E,X)=0 ~\forall E\in \cE. \]

One can show that $\cE^\bot$ is also triangulated and the functor
defined as the composition 
\[ \cE^\bot\ra \cC \ra \cC/\cE \] 
is a triangulated equivalence. One usually denotes this situation by
$\cC=\langle\cE^\bot,\cE\rangle$ to illustrate the fact that every
object in $\cC$ can be obtained from $\cE$ and $\cE^\bot$ by taking
cones and shifts.
\medskip

\noindent
We can make a similar definition of the notion of being admissible in
the dg setting.

\begin{Definition} 
Let $\cD$ be a pre-triangulated dg category and let $\cE$ be a full
subcategory of $\cD$ which is also pre-triangulated. Then $\cE$ is
said to be right (left) admissible in $\cD$ if $H^0(\cE)$ is
right (left) admissible in $H^0(\cD)$. 

We shall denote by $\cE^\bot$ the full subcategory of $\cD$ consists
of those objects in $\cD$ that are inside $H^0(\cE)^\bot$, or,
equivalently, $\cE^\bot$ consists of those objects $X$ such that
$\Hom_\cD(E,X)$ is acyclic for any $E\in \cE$.  \medskip

\noindent
As before  we write
\[ \cD=<\cE^\bot,\cE>.\] 
\end{Definition}
\medskip

\noindent
We end this section with an easy Lemma that is used in the proofs in
Section~\ref{sec:dg}.

\begin{Lemma}
\label{lem:adm} Let $\cE$ be a right (left) admissible pre-triangulated
subcategory of the pre-triangulated category $\cD$. Then
$\cE^\bot$ (${}^\bot\cE$) is quasi-equivalent to the dg quotient
$\cD/\cE$.
\end{Lemma}

\begin{Proof} 
Consider the composition of morphisms in $\Hqe$
 \[\cE^\bot  \hookrightarrow \cD \rightarrow \cD/\cE. \] 
To show that this map is an isomorphism in $\Hqe$ we only need to check the statement in the homotopy category. Consider the induced map on the
corresponding homotopy categories 
\[H^0(\cE^\bot) \hookrightarrow H^0(\cD) \rightarrow H^0(\cD/\cE).\] 
Since both $\cE$ and $\cD$ are pre-triangulated we have 
\[ H^0(\cD/\cE) \cong (\cD/\cE)^\tr.\]
By Theorem~\ref{thm:quo} we have 
\[ (\cD/\cE)^\tr\cong \cD^\tr/\cE^\tr,\] 
and the fact that the morphism
\[ H^0(\cE^\bot)\rightarrow \cD^\tr/\cE^\tr \] 
is an isomorphism is the known result for triangulated categories.
\qed
\end{Proof}

\end{document}